\documentclass[11pt]{amsart}
\usepackage[utf8]{inputenc}
\usepackage[T1]{fontenc}
\usepackage{bm,amsmath,amsthm,amssymb}
\usepackage{geometry}
\usepackage{graphicx}
\usepackage{tikz}
\usetikzlibrary{patterns}
\usepackage{amsfonts}
\usepackage{array}
\usepackage{enumerate}
\usepackage{enumitem}
\usepackage{float}
\usepackage{lscape}
\usepackage[english]{babel}
\usepackage[all]{xy}
\usepackage{authblk}
\usepackage{hyperref}

\theoremstyle{plain}
\newtheorem{fact}{Fact}[section]
\newtheorem{theo}[fact]{Theorem}
\newtheorem{lem}[fact]{Lemma}
\newtheorem{defi}[fact]{Definition}

\newtheorem{prop}[fact]{Proposition}
\newtheorem{rmk}[fact]{Remark}

\title{A construction of the generalized higher cluster category arising from an $(m+2)$-angulation of a marked surface}

\begin{document}

\maketitle

\begin{center}
\Large	\textit{Lucie Jacquet-Malo \\
	lucie.jacquet.malo@u-picardie.fr}
\end{center}

\begin{abstract}
	In this article, we study the $(m+2)$-angulations on a Riemann surface, characterized with its boundary components, punctures, and gender. We count the number of arcs in such a surface, and associate a graded quiver with superpotential associated with an $(m2)$-angulation. We show the compatibility between the flip of an $(m+2)$-angulation and the flip in the unpunctured case.
\end{abstract}



\section*{Introduction}

In 2002, Fomin and Zelevinsky introduced cluster algebras in \cite{FZ} to establish a combinatorial framework for studying canonical bases. Since then, cluster algebras have been found to be deeply connected to various mathematical areas, including Calabi-Yau algebras, integrable systems, Poisson geometry, and quiver representations. To further develop this concept, Buan, Marsh, Reineke, Reiten, and Todorov introduced cluster categories in \cite{BMRRT}, with Caldero, Chapoton, and Schiffler providing a specific case for $A_n$ in \cite{CCS}. These developments enabled the categorification of mutations in cluster algebras through tilting theory. For a gentle introduction to cluster categories, refer to Keller's article \cite{Kel01}.

\vspace{20pt}

The cluster category is defined as follows: Let $K$ be a field, $Q$ an acyclic quiver, and $\mathcal{D}^b(KQ)$ the bounded derived category of $Q$. The cluster category is the orbit of $\mathcal{D}^b(KQ)$ under the functor $\tau^{-1}[1]$, where $\tau$ denotes the Auslander-Reiten translation, and $[1]$ the shift functor. Keller demonstrated in \cite{Kel03} that the cluster category is triangulated, with the shift functor $[1]$. He also proved that for sufficiently well-behaved endofunctors, the orbit category of a derived category is triangulated. This result led to the concept of higher cluster categories, defined as 
\[\mathcal{D}^b(KQ)/\tau^{-1}[m].\] Thomas in \cite{Tho} formally defined higher cluster categories, showing they serve a similar role as cluster categories but concerning $m$-clusters (as defined by Fomin and Reading in \cite{FR}). Wraalsen and Zhou/Zhu in \cite{W} and \cite{ZZ} exrtended many properties of cluster categories to higher cluster categories. For instance, they demonstrated that any $m$-rigid object $X$ with $n-1$ nonisomorphic indecomposable summands has exactly $m+1$ complements, meaning nonisomorphic indecomposable objects $Y$ such that $X \oplus Y$ forms an $m$-cluster-tilting object.

For certain classes of quivers, it is possible to construct geometric realizations of (higher) cluster categories. Caldero, Chapoton, and Schiffler achieved this for $A_n$ in cluster categories, while Schiffler extended it to $D_n$ in \cite{Sch}. Baur and Marsh generalized these results to higher cluster categories in \cite{BM01} and \cite{BM02}. In these cases, the Auslander-Reiten quiver of the higher cluster category can be realized as a connected component of a geometrically constructed category.

However, this geometric realization is not possible in Euclidean cases, which are representation-infinite. In these cases, the Auslander-Reiten quiver of the higher cluster category is infinite and consists of three repeating main parts. Torkildsen explored the case of $\tilde{A}$ in \cite{Tor}, and Baur and Torkildsen provided a complete geometric realization of this case in \cite{BauTor}. Case $\tilde{D}$ has been treated in \cite{LJM} and \cite{LJM2}.

In a parallel way, Amiot in \cite{A} introduced higher cluster categories, which have been generalized by Guo in \cite{G} for higher cluster categories (which are $(m+1)$-Calabi-Yau). Note that the gerenalized higher cluster categories arise from a graded quiver with superpotential.

The paper is aimed to give a geometric realization of $m$-cluster categories for all types of surfaces. Unfortunately, the usual descriptions do not fit perfectly. To be more accurate, there is no compatibility between the flip of an $(m+2)$-angulation and the flip of a quiver with superpotential, in the punctured case.

The paper is organized as follows :

In section $1$, we recall some important notions on higher cluster categories, mutation of $m$-rigid objects and quivers with superpotential.

Section $2$ is devoted to introduce the $(m+2)$-angulation associated with a Riemann surface. The main Theorem of this section is \ref{th:arcs}, which counts the number of arcs in an $(m+2)$-angulation.

In section $3$, we associate a graded quiver with super potential with an $(m+2)$-angulation, and show the compatibility with the flip.

Finally, in section $4$, we build a higher cluster category from the quiver given in the previous section, using the results of Keller and Guo in \cite{Keldefo} and \cite{G}.

\section{Preliminaries}

Notations:

Throughout this paper, we fix a field K and an acyclic finite quiver $Q$.
In the remaining of the paper, $n$ and $m$ are integers, where $n$ is the number of vertices of $Q$, $n \geq 4$. We note that all the results apply to the cases $A_1$, $A_2$, and $A_3$ using exactly the same arguments.

If $A$ is an object in a category ${\mathcal{C}}$, $A^\perp$ is the class of all objects $X$ such that \[{\mathrm{Ext}}_\mathcal{C}^i(A,X)=0 \text{ for all } i \in \{1,\cdots,m \}.\]

The category ${\mathrm{mod}}(KQ)$ is the category of finitely generated right modules over the path algebra $KQ$. The letter $\tau$ stands for the Auslander-Reiten translation. We write $[1]$ for the shift functor in the bounded derived category $\mathcal{D}^b(KQ)$. For any further information about representation theory of associative algebras, see the book written by Assem, Simson and Skowronski, \cite{ASS}.
The symbol $\check{M}$ replaces $\Sigma^n R \mathrm{Hom}_{A^e}(M,A^e)$, where $A^e=A \otimes A^{\mathrm{op}}$. The hooks $[,]$ denote the supercommutator. The couple $(S,M)$ is a surface $S$ with a set of marked points $M$ and eventually punctures. When we talk about a quotient, the application $\pi$ denotes the canonical projection.

:\subsection{Higher cluster categories}

In 2006, in order to understand better the notion of cluster algebras, Buan, Marsh, Reineke, Reiten and Todorov in \cite{BMRRT} defined the cluster category of an acyclic quiver in the following way:

If $Q$ is an acyclic quiver, let ${\mathcal{D}}^b(KQ)$ be the bounded derived category of the category ${\mathrm{mod}}~KQ$. The category ${\mathcal{C}}_Q$ is the orbit category in the sense of the derived category under the functor $\tau^{-1}[1]$ in the sense of Keller in \cite{Kel03}.

Cluster categories give a real interpretation of clusters in a cluster algebra in terms of tilting objects. To be precise, the cluster variables of the algebra are in $1-1$ correspondence with the indecomposables objects in ${\mathcal{C}}_Q$, and the clusters are in $1-1$ correspondence with the basic tilting objects in ${\mathcal{C}}_Q$.

It is known from \cite{BMRRT} that ${\mathcal{C}}_Q$ is Krull-Schmidt. Since $\tau=[1]$, we have that ${\mathcal{C}}_Q$ is $2$-Calabi-Yau, and Keller in \cite{Kel03} has shown that it was a triangulated category.

If $m$ is a nonzero integer, we can also define the higher cluster category
\[ {\mathcal{C}}^m_Q={\mathcal{D}}^b(KQ)/\tau^{-1}[m]. \]
Again, the higher cluster category is Krull-Schmidt, $(m+1)$-Calabi-Yau, and triangulated.

\begin{defi}
Let $T$ be an object in the category ${\mathcal{C}}^m_Q$. Then $T$ is $m$-rigid if
\[ {\mathrm{Ext}}^i(T,T)=0~\forall i \in \{1,\cdots,m \}. \]
\end{defi}

\begin{defi}\cite{KR}
Let $T$ be an object in the category ${\mathcal{C}}^m_Q$. Then $T$ is $m$-cluster-tilting when we have the following equivalence: \[ X \text{ is in } {\mathrm{add}}~T \iff {\mathrm{Ext}}^i_{{\mathcal{C}}_Q}(T,X)=0~\forall i \in \{ 1,\cdots,m \}. \]
\end{defi}

It is known from Zhu in \cite{Z}, that $T$ is an $m$-cluster-tilting object if and only if $T$ has $n$ indecomposable direct summands (up to isomorphism) and is $m$-rigid. So, for $T=\bigoplus T_i$ an $m$-cluster-tilting object, where each $T_i$ is indecomposable, let us define the almost complete $m$-rigid object ${\overline{T}}=T/T_k$ (where $T_k$ is an indecomposable summand of $T$). There are, up to isomorphism, $m+1$ complements of the almost $m$-cluster-tilting object $\overline{T}$, denoted by $T_k^{(c)}$, for $c \in \{ 0,\cdots, m \}$. Iyama and Yoshino in \cite{IY} showed in 2008 the following theorem:

\begin{theo}\cite{IY}
There are $m+1$ exchange triangles:
\[
\xymatrix{
T_k^{(c)} \ar^{f_k^{(c)}}[r] & B_k^{(c)} \ar^{g_k^{(c+1)}}[r] & T_k^{(c+1)} \ar^{h_k^{(c+1)}}[r] & T_k^{(c)}[1] }
\]
Here, the $B_k^{(c)}$ are in ${\mathrm{add}} \overline{T}$, the maps $f_k^{(c)}$ (respectively $g_k^{(c+1)}$) are minimal left (respectively right) ${\mathrm{add}} \overline{T}$-approximations, hence, not split monomorphism or split epimorphism.
\end{theo}

\subsection{Generalized higher cluster categories (\cite{G})}

We are now interested in generalized higher cluster categories introduced by Guo in \cite{G}, in the way Amiot defined generalized cluster categories in \cite{A}.

\begin{defi}\cite{Gin}
Let $A$ be an algebra over $K$. Then $A$ is homologically smooth if, as a bimodule, it admits a finite resolution by finitely generated projective bimodules.

Moreover, $A$ is $n$-Calabi-Yau as a bimodule if it is homologically smooth and, in the derived category of $A$-bimodules, we have an isomorphism:

\[ f: \Sigma^n R \mathrm{Hom}_{A^e}(A,A^e) \to A \]

such that $\Sigma^n R \mathrm{Hom}_{A^e}(f,A^e)=f$.
\end{defi}

In her article, Guo showed the following theorem:

\begin{theo}\label{th:guo}\cite{G}
Let $A$ be a differential graded (=dg) algebra having the following properties:
\begin{enumerate}
\item The algebra $A$ is homologically smooth,
\item for all $p > 0$ we have $H_p(A)=0$,
\item the zeroth homology $H_0(A)$ is finite dimensional,
\item the algebra $A$ is $(m+2)$-Calabi-Yau as a bimodule.
\end{enumerate}
Then
\begin{enumerate}
\item The category $C_A=\mathrm{per}A/\mathcal{D}^b(A)$ is Hom-finite and $(m+1)$-Calabi-Yau,
\item the object $T=\pi A$ is an $m$-cluster-tilting object,
\item we have an isomorphism $\mathrm{End}(T) \simeq H_0(A)$.
\end{enumerate}
\end{theo}

In the case where $A=KQ$ for an acyclic quiver $Q$, we find back the classical higher cluster category.

\subsection{Ginzburg dg categories}
Let $Q$ be a graded $k$-quiver, whose set of objects is finite and $Q(x,y)$ is a finitely generated graded projective $k$-module for any objects $x$ and $y$. Let $W$ be a superpotential of degree $n-3$, it means roughly a linear combination of cycles considered up to cyclic permutation with signs. Let $\mathcal{R}$ be the discrete category on $Q_0$.

\begin{defi}\cite{Keldefo}
The Ginzburg dg category $\Gamma_n(Q,W)$ is defined as the tensor category over $\mathcal{R}$ of the bimodule
\[ \tilde{Q}=Q \oplus \check{Q}[n-2] \oplus \mathcal{R}[n-1] \]
endowed with the unique differential which

\begin{enumerate}
\item vanishes on $Q$,
\item takes the element $v_i^*$ of $\check{Q}[n-2]$ to the cyclic derivative $\partial_{v_i}W$
\item takes the element $\mathrm{id}_x$ of $\mathcal{R}[n-1]$ to $(-1)^n \mathrm{id}_x (\sum [v_i,v_i^*])\mathrm{id}_x$
\end{enumerate}
\end{defi}

This definition leads to the following theorem:

\begin{theo}\cite{Keldefo}\label{th:kel}
The dg Ginzburg category is homologically smooth and $n$-Calabi-Yau.
\end{theo}

\subsection{Quivers with superpotential and their mutations}
Let us define now a potential on an algebra $A$.

\begin{defi}\cite{DWZ}
Let $A$ be an algebra. Let $R\langle \langle A \rangle \rangle=\prod_{d=1}^\infty A^d$. Then a potential $W$ on $Q$ is an element of the closed vector space
\[ R\langle \langle A \rangle \rangle_{cyc}=\prod_{d=1}^\infty A^d_{cyc}, \]
where $A^d_{cyc}$ is the algebra generated by the cycles of order $d$.

For any $\xi$ in $A^\star$, where $A^\star$ is the dual bimodule of $A$, we define the cyclic derivative $\partial_\xi$ as the continuous $K$-linear map $R\langle \langle A \rangle \rangle_{cyc} \to R\langle \langle A \rangle \rangle$, acting on paths by
\[ \partial_\xi(a_1,\cdots,a_d)=\sum_{k=1}^d \xi(a_k) a_{k+1}\cdots a_d a_1\cdots a_{k-1} \]

For any potential $W$, we define its Jacobian ideal $J(W)$ as the closure of the ideal in $R\langle \langle A \rangle \rangle$ generated by the elements $\partial_\xi(W)$. We call by Jacobian algebra of $W$ the quotient $R\langle \langle A \rangle \rangle/J(W)$.
\end{defi}

\begin{defi}
Two potentials $W$ and $W'$ are called cyclically equivalent if $W-W'$ lies in the closure of the span of all elements of the form $a_1\cdots a_d -a_2\cdots a_d a_1$, where $a_1\cdots a_d$ is a cycle.
\end{defi}

We now define a quiver with potential.

\begin{defi}
Let $Q$ be a quiver and let $A$ be its path algebra. Then the couple $(Q,W)$ is a quiver with potential if it satisfies the following conditions:
\begin{enumerate}
\item The quiver $Q$ has no loops,
\item no two cyclically equivalent cyclic paths appear in $W$.
\end{enumerate}
\end{defi}

Derksen, Weyman and Zelevinsky introduced the mutation of a quiver with potential. We report to the article \cite[Section 5]{DWZ} for a complete definition.

\section{Definition of $(m+2)$-angulations of a surface}

\subsection{Number of arcs in an $(m+2)$-angulation}

We first define $(m+2)$-angulations, in order to find the number of arcs belonging to an $(m+2)$-angulation in function of the characteristics of the surface.

\begin{defi}
An arc in $(S,M)$ is the homotopy class relative to endpoints of a path between two marked points, or two punctures, or one marked point and one puncture.

Two arcs do not cross if their homotopy classes contain representatives which do not cross.

Arcs around a puncture whose endpoints coincide are called loops.
\end{defi}

By an arc, we mean a representative in the class of homotopy of its endpoints.

\begin{defi}
Let $(S,M)$ be a marked surface with punctures. Then an $(m+2)$-angulation is a set of noncrossing arcs cutting the universal cover of $(S,M)$ into pseudo-$(m+2)$-angles, it means figures delimited by $m+2$ paths.
\end{defi}

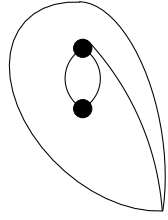
\begin{figure}[!h]
\centering
\begin{tikzpicture}[scale=0.4]
\draw [fill=black,fill opacity=1.0] (2,2) circle (0.3cm);
\draw [fill=black,fill opacity=1.0] (2,0) circle (0.3cm);
\draw [shift={(1.43,1)}] plot[domain=-1.06:1.06,variable=\t]({1*1.15*cos(\t r)+0*1.15*sin(\t r)},{0*1.15*cos(\t r)+1*1.15*sin(\t r)});
\draw [shift={(2.57,1)}] plot[domain=2.09:4.2,variable=\t]({1*1.15*cos(\t r)+0*1.15*sin(\t r)},{0*1.15*cos(\t r)+1*1.15*sin(\t r)});
\draw [shift={(3.11,0.07)}] plot[domain=-0.66:2.58,variable=\t]({-0.72*4.13*cos(\t r)+-0.69*2.74*sin(\t r)},{0.69*4.13*cos(\t r)+-0.72*2.74*sin(\t r)});
\draw [shift={(2.45,-1.07)}] plot[domain=4.11:6.2,variable=\t]({-0.16*4.65*cos(\t r)+-0.99*2.16*sin(\t r)},{0.99*4.65*cos(\t r)+-0.16*2.16*sin(\t r)});
\draw [shift={(-4.66,-4.06)}] plot[domain=0.07:0.74,variable=\t]({1*9.3*cos(\t r)+0*9.3*sin(\t r)},{0*9.3*cos(\t r)+1*9.3*sin(\t r)});
\end{tikzpicture}
\caption{Example of an $(m+2)$-gon}
\label{fig:oeil}
\end{figure}

An admissible arc, or an $m$-diagonal, is an arc that appears in an $(m+2)$-angulation.

\begin{figure}[!h]
\centering
\begin{tikzpicture}[scale=0.6]
\draw [rotate around={0:(4.73,0)}] (4.73,0) ellipse (5.94cm and 3.6cm);
\draw [shift={(4.57,0.94)}] plot[domain=0.34:3.12,variable=\t]({-1*3.79*cos(\t r)+-0.09*2.16*sin(\t r)},{0.09*3.79*cos(\t r)+-1*2.16*sin(\t r)});
\draw [shift={(4.74,-0.07)}] plot[domain=0.11:3.08,variable=\t]({1*3.43*cos(\t r)+0.01*1.23*sin(\t r)},{-0.01*3.43*cos(\t r)+1*1.23*sin(\t r)});
\draw [fill=black,pattern=north east lines] (4.16,-2.34) circle (0.3cm);
\draw [shift={(5,-1.95)},color=red]  plot[domain=2.09:3.25,variable=\t]({1*0.84*cos(\t r)+0*0.84*sin(\t r)},{0*0.84*cos(\t r)+1*0.84*sin(\t r)});
\draw [shift={(6.07,-2.38)},color=red]  plot[domain=3.28:3.83,variable=\t]({1*1.93*cos(\t r)+0*1.93*sin(\t r)},{0*1.93*cos(\t r)+1*1.93*sin(\t r)});
\draw [shift={(2.03,-2.41)},dash pattern=on 3pt off 3pt,color=red]  plot[domain=-0.44:0.44,variable=\t]({1*2.81*cos(\t r)+0*2.81*sin(\t r)},{0*2.81*cos(\t r)+1*2.81*sin(\t r)});
\draw [color=green] (4.16,-2.04)-- (6.29,-1.09);
\draw [dash pattern=on 5pt off 5pt,color=green] (6.29,-1.09)-- (8.43,-2.82);
\draw [shift={(4.19,-0.23)},color=blue]  plot[domain=1.53:5.17,variable=\t]({-1*4.94*cos(\t r)+0.07*2.4*sin(\t r)},{-0.07*4.94*cos(\t r)+-1*2.4*sin(\t r)});
\draw [shift={(4.51,0.04)},color=blue]  plot[domain=-0.91:1.51,variable=\t]({-1*4.4*cos(\t r)+-0.03*2.09*sin(\t r)},{0.03*4.4*cos(\t r)+-1*2.09*sin(\t r)});
\draw [shift={(4.5,-0.48)},color=green]  plot[domain=2.32:5.6,variable=\t]({-0.99*5.38*cos(\t r)+0.11*3.74*sin(\t r)},{-0.11*5.38*cos(\t r)+-0.99*3.74*sin(\t r)});
\draw [shift={(4.47,0.78)},color=green]  plot[domain=-0.66:1.38,variable=\t]({-0.97*5.2*cos(\t r)+0.24*2.71*sin(\t r)},{-0.24*5.2*cos(\t r)+-0.97*2.71*sin(\t r)});
\begin{scriptsize}
\draw [color=black] (4.16,-2.04)-- ++(-1.5pt,-1.5pt) -- ++(3.0pt,3.0pt) ++(-3.0pt,0) -- ++(3.0pt,-3.0pt);
\draw [color=black] (4.16,-2.64)-- ++(-1.5pt,-1.5pt) -- ++(3.0pt,3.0pt) ++(-3.0pt,0) -- ++(3.0pt,-3.0pt);
\end{scriptsize}
\end{tikzpicture}
$\to$
\begin{tikzpicture}[scale=0.8]
\draw (-3,4)-- (-3,-5);
\draw (0,4)-- (0,-5);
\draw (3,-5)-- (3,4);
\draw (-3,4)-- (3,4);
\draw (-3,-5)-- (3,-5);
\draw (-3,1)-- (3,1);
\draw (-3,-2)-- (3,-2);
\draw [fill=black,pattern=north east lines] (-1.5,-0.5) circle (0.3cm);
\draw [fill=black,pattern=north east lines] (-1.5,-3.5) circle (0.3cm);
\draw [fill=black,pattern=north east lines] (-1.5,2.5) circle (0.3cm);
\draw [fill=black,pattern=north east lines] (1.5,2.5) circle (0.3cm);
\draw [fill=black,pattern=north east lines] (1.5,-0.5) circle (0.3cm);
\draw [fill=black,pattern=north east lines] (1.5,-3.5) circle (0.3cm);
\draw [shift={(-0.97,0.45)},color=blue]  plot[domain=4.31:5.51,variable=\t]({1*1.36*cos(\t r)+0*1.36*sin(\t r)},{0*1.36*cos(\t r)+1*1.36*sin(\t r)});
\draw [shift={(0.97,-1.45)},color=blue]  plot[domain=1.17:2.37,variable=\t]({1*1.36*cos(\t r)+0*1.36*sin(\t r)},{0*1.36*cos(\t r)+1*1.36*sin(\t r)});
\draw [shift={(-0.97,3.45)},color=blue]  plot[domain=4.31:5.51,variable=\t]({1*1.36*cos(\t r)+0*1.36*sin(\t r)},{0*1.36*cos(\t r)+1*1.36*sin(\t r)});
\draw [shift={(-0.97,-2.55)},color=blue]  plot[domain=4.31:5.51,variable=\t]({1*1.36*cos(\t r)+0*1.36*sin(\t r)},{0*1.36*cos(\t r)+1*1.36*sin(\t r)});
\draw [shift={(0.97,1.55)},color=blue]  plot[domain=1.17:2.37,variable=\t]({1*1.36*cos(\t r)+0*1.36*sin(\t r)},{0*1.36*cos(\t r)+1*1.36*sin(\t r)});
\draw [shift={(0.97,-4.45)},color=blue]  plot[domain=1.17:2.37,variable=\t]({1*1.36*cos(\t r)+0*1.36*sin(\t r)},{0*1.36*cos(\t r)+1*1.36*sin(\t r)});
\draw [shift={(2.03,3.45)},color=blue]  plot[domain=4.31:5.51,variable=\t]({1*1.36*cos(\t r)+0*1.36*sin(\t r)},{0*1.36*cos(\t r)+1*1.36*sin(\t r)});
\draw [shift={(2.03,0.45)},color=blue]  plot[domain=4.31:5.51,variable=\t]({1*1.36*cos(\t r)+0*1.36*sin(\t r)},{0*1.36*cos(\t r)+1*1.36*sin(\t r)});
\draw [shift={(2.03,-2.55)},color=blue]  plot[domain=4.31:5.51,variable=\t]({1*1.36*cos(\t r)+0*1.36*sin(\t r)},{0*1.36*cos(\t r)+1*1.36*sin(\t r)});
\draw [shift={(-2.03,-4.45)},color=blue]  plot[domain=1.17:2.37,variable=\t]({1*1.36*cos(\t r)+0*1.36*sin(\t r)},{0*1.36*cos(\t r)+1*1.36*sin(\t r)});
\draw [shift={(-2.03,-1.45)},color=blue]  plot[domain=1.17:2.37,variable=\t]({1*1.36*cos(\t r)+0*1.36*sin(\t r)},{0*1.36*cos(\t r)+1*1.36*sin(\t r)});
\draw [shift={(-2.03,1.55)},color=blue]  plot[domain=1.17:2.37,variable=\t]({1*1.36*cos(\t r)+0*1.36*sin(\t r)},{0*1.36*cos(\t r)+1*1.36*sin(\t r)});
\draw [color=red] (-1.5,2.8)-- (-1.5,4);
\draw [color=red] (-1.5,2.2)-- (-1.5,-0.2);
\draw [color=red] (-1.5,-0.8)-- (-1.5,-3.2);
\draw [color=red] (-1.5,-3.8)-- (-1.5,-5);
\draw [color=red] (1.5,2.8)-- (1.5,4);
\draw [color=red] (1.5,2.2)-- (1.5,-0.2);
\draw [color=red] (1.5,-0.8)-- (1.5,-3.2);
\draw [color=red] (1.5,-3.8)-- (1.5,-5);
\draw [shift={(-1.78,0.96)},color=green]  plot[domain=-1.41:0.02,variable=\t]({1*1.78*cos(\t r)+0*1.78*sin(\t r)},{0*1.78*cos(\t r)+1*1.78*sin(\t r)});
\draw [shift={(1.78,1.04)},color=green]  plot[domain=1.73:3.16,variable=\t]({1*1.78*cos(\t r)+0*1.78*sin(\t r)},{0*1.78*cos(\t r)+1*1.78*sin(\t r)});
\draw [shift={(-1.78,3.96)},color=green]  plot[domain=-1.41:0.02,variable=\t]({1*1.78*cos(\t r)+0*1.78*sin(\t r)},{0*1.78*cos(\t r)+1*1.78*sin(\t r)});
\draw [shift={(-1.78,-2.04)},color=green]  plot[domain=-1.41:0.02,variable=\t]({1*1.78*cos(\t r)+0*1.78*sin(\t r)},{0*1.78*cos(\t r)+1*1.78*sin(\t r)});
\draw [shift={(1.78,-1.96)},color=green]  plot[domain=1.73:3.16,variable=\t]({1*1.78*cos(\t r)+0*1.78*sin(\t r)},{0*1.78*cos(\t r)+1*1.78*sin(\t r)});
\draw [shift={(1.78,-4.96)},color=green]  plot[domain=1.73:3.16,variable=\t]({1*1.78*cos(\t r)+0*1.78*sin(\t r)},{0*1.78*cos(\t r)+1*1.78*sin(\t r)});
\draw [shift={(-1.22,1.04)},color=green]  plot[domain=1.73:3.16,variable=\t]({1*1.78*cos(\t r)+0*1.78*sin(\t r)},{0*1.78*cos(\t r)+1*1.78*sin(\t r)});
\draw [shift={(1.22,3.96)},color=green]  plot[domain=-1.41:0.02,variable=\t]({1*1.78*cos(\t r)+0*1.78*sin(\t r)},{0*1.78*cos(\t r)+1*1.78*sin(\t r)});
\draw [shift={(1.22,0.96)},color=green]  plot[domain=-1.41:0.02,variable=\t]({1*1.78*cos(\t r)+0*1.78*sin(\t r)},{0*1.78*cos(\t r)+1*1.78*sin(\t r)});
\draw [shift={(-1.22,-1.96)},color=green]  plot[domain=1.73:3.16,variable=\t]({1*1.78*cos(\t r)+0*1.78*sin(\t r)},{0*1.78*cos(\t r)+1*1.78*sin(\t r)});
\draw [shift={(1.22,-2.04)},color=green]  plot[domain=-1.41:0.02,variable=\t]({1*1.78*cos(\t r)+0*1.78*sin(\t r)},{0*1.78*cos(\t r)+1*1.78*sin(\t r)});
\draw [shift={(-1.22,-4.96)},color=green]  plot[domain=1.73:3.16,variable=\t]({1*1.78*cos(\t r)+0*1.78*sin(\t r)},{0*1.78*cos(\t r)+1*1.78*sin(\t r)});
\begin{scriptsize}
\draw [color=black] (-1.5,2.8)-- ++(-1.5pt,-1.5pt) -- ++(3.0pt,3.0pt) ++(-3.0pt,0) -- ++(3.0pt,-3.0pt);
\draw [color=black] (-1.5,2.2)-- ++(-1.5pt,-1.5pt) -- ++(3.0pt,3.0pt) ++(-3.0pt,0) -- ++(3.0pt,-3.0pt);
\draw [color=black] (-1.5,-0.2)-- ++(-1.5pt,-1.5pt) -- ++(3.0pt,3.0pt) ++(-3.0pt,0) -- ++(3.0pt,-3.0pt);
\draw [color=black] (-1.5,-0.8)-- ++(-1.5pt,-1.5pt) -- ++(3.0pt,3.0pt) ++(-3.0pt,0) -- ++(3.0pt,-3.0pt);
\draw [color=black] (-1.5,-3.2)-- ++(-1.5pt,-1.5pt) -- ++(3.0pt,3.0pt) ++(-3.0pt,0) -- ++(3.0pt,-3.0pt);
\draw [color=black] (-1.5,-3.8)-- ++(-1.5pt,-1.5pt) -- ++(3.0pt,3.0pt) ++(-3.0pt,0) -- ++(3.0pt,-3.0pt);
\draw [color=black] (1.5,2.8)-- ++(-1.5pt,-1.5pt) -- ++(3.0pt,3.0pt) ++(-3.0pt,0) -- ++(3.0pt,-3.0pt);
\draw [color=black] (1.5,2.2)-- ++(-1.5pt,-1.5pt) -- ++(3.0pt,3.0pt) ++(-3.0pt,0) -- ++(3.0pt,-3.0pt);
\draw [color=black] (1.5,-0.8)-- ++(-1.5pt,-1.5pt) -- ++(3.0pt,3.0pt) ++(-3.0pt,0) -- ++(3.0pt,-3.0pt);
\draw [color=black] (1.5,-3.8)-- ++(-1.5pt,-1.5pt) -- ++(3.0pt,3.0pt) ++(-3.0pt,0) -- ++(3.0pt,-3.0pt);
\draw [color=black] (1.5,-0.2)-- ++(-1.5pt,-1.5pt) -- ++(3.0pt,3.0pt) ++(-3.0pt,0) -- ++(3.0pt,-3.0pt);
\draw [color=black] (1.5,-3.2)-- ++(-1.5pt,-1.5pt) -- ++(3.0pt,3.0pt) ++(-3.0pt,0) -- ++(3.0pt,-3.0pt);
\end{scriptsize}
\end{tikzpicture}
\caption{A $4$-angulation of the torus with one boundary component and two marked points on it. On the right side we can see the universal cover of it, and the pseudo $(m+2)$-gons.}
\label{fig:ang}
\end{figure}
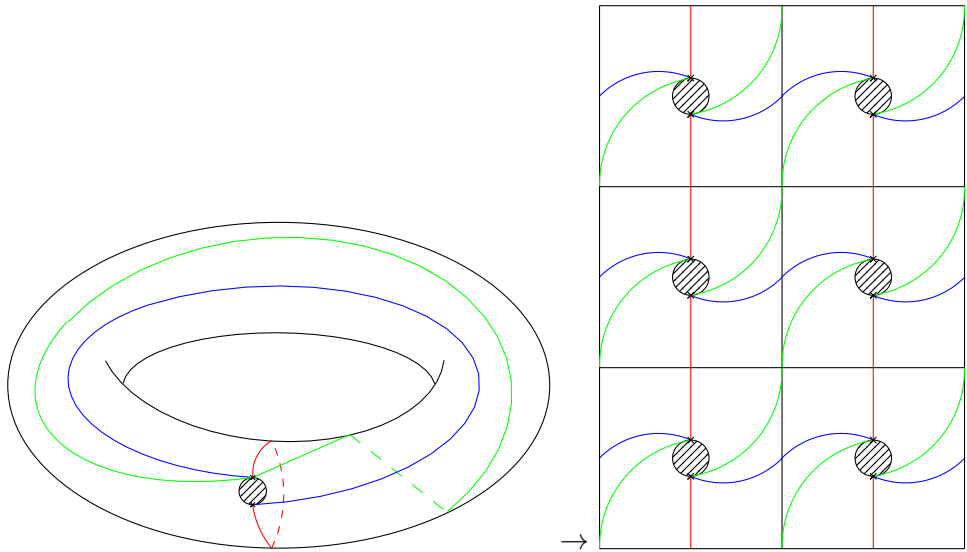

We illustrate the notion of $(m+2)$-angulation on a torus on figure \ref{fig:ang}.

\begin{figure}[!h]
\centering
\begin{tikzpicture}[scale=0.8]
\draw [fill=black,fill opacity=1.0] (6.3,-2.23) circle (0.2cm);
\draw [rotate around={0:(4.73,0)}] (4.73,0) ellipse (5.94cm and 3.6cm);
\draw [shift={(4.57,0.94)}] plot[domain=0.34:3.12,variable=\t]({-1*3.79*cos(\t r)+-0.09*2.16*sin(\t r)},{0.09*3.79*cos(\t r)+-1*2.16*sin(\t r)});
\draw [shift={(4.74,-0.07)}] plot[domain=0.11:3.08,variable=\t]({1*3.43*cos(\t r)+0.01*1.23*sin(\t r)},{-0.01*3.43*cos(\t r)+1*1.23*sin(\t r)});
\draw [fill=black,pattern=north east lines] (4.16,-2.34) circle (0.3cm);
\draw [shift={(5.84,-4.99)},color=violet]  plot[domain=1.4:2.05,variable=\t]({1*3*cos(\t r)+0*3*sin(\t r)},{0*3*cos(\t r)+1*3*sin(\t r)});
\draw [shift={(5.6,1.47)},color=orange]  plot[domain=4.42:4.9,variable=\t]({1*3.96*cos(\t r)+0*3.96*sin(\t r)},{0*3.96*cos(\t r)+1*3.96*sin(\t r)});
\draw [shift={(7.38,-2.41)},color=red]  plot[domain=2.72:3.56,variable=\t]({1*2.92*cos(\t r)+0*2.92*sin(\t r)},{0*2.92*cos(\t r)+1*2.92*sin(\t r)});
\draw [shift={(1.65,-2.41)},dash pattern=on 5pt off 5pt,color=red]  plot[domain=-0.37:0.37,variable=\t]({1*3.28*cos(\t r)+0*3.28*sin(\t r)},{0*3.28*cos(\t r)+1*3.28*sin(\t r)});
\draw [color=blue] (4.46,-2.33)-- (7.08,-3.31);
\draw [dash pattern=on 5pt off 5pt,color=blue] (7.08,-3.31)-- (7.66,-0.49);
\draw [shift={(4.69,0.59)},color=blue]  plot[domain=2.28:6.09,variable=\t]({-1*4.53*cos(\t r)+0.01*1.46*sin(\t r)},{-0.01*4.53*cos(\t r)+-1*1.46*sin(\t r)});
\draw [shift={(3.44,-0.07)},color=blue]  plot[domain=-0.27:1.8,variable=\t]({-0.99*3.43*cos(\t r)+-0.16*1.49*sin(\t r)},{0.16*3.43*cos(\t r)+-0.99*1.49*sin(\t r)});
\draw [shift={(3.9,-2.2)},color=blue]  plot[domain=-0.22:1.46,variable=\t]({1*0.57*cos(\t r)+0*0.57*sin(\t r)},{0*0.57*cos(\t r)+1*0.57*sin(\t r)});
\draw [shift={(6.09,-3.73)},color=green]  plot[domain=1.5:2.43,variable=\t]({1*2.15*cos(\t r)+0*2.15*sin(\t r)},{0*2.15*cos(\t r)+1*2.15*sin(\t r)});
\draw [shift={(5.19,0.57)},color=green]  plot[domain=1.81:5.08,variable=\t]({-1*4.7*cos(\t r)+-0.05*2.18*sin(\t r)},{0.05*4.7*cos(\t r)+-1*2.18*sin(\t r)});
\draw [shift={(3.7,0.12)},color=green]  plot[domain=-1.63:0.39,variable=\t]({-0.99*4.57*cos(\t r)+0.15*2.55*sin(\t r)},{-0.15*4.57*cos(\t r)+-0.99*2.55*sin(\t r)});
\draw [shift={(3.46,2.56)},color=green]  plot[domain=3.96:4.78,variable=\t]({1*5.52*cos(\t r)+0*5.52*sin(\t r)},{0*5.52*cos(\t r)+1*5.52*sin(\t r)});
\draw [shift={(3.86,-2.35)},color=green]  plot[domain=-1.65:0.04,variable=\t]({1*0.61*cos(\t r)+0*0.61*sin(\t r)},{0*0.61*cos(\t r)+1*0.61*sin(\t r)});
\begin{scriptsize}
\draw [color=black] (4.46,-2.33)-- ++(-1.5pt,-1.5pt) -- ++(3.0pt,3.0pt) ++(-3.0pt,0) -- ++(3.0pt,-3.0pt);
\end{scriptsize}
\end{tikzpicture}
\caption{Here, g=1, b=1, c=1, p=1, so five arcs compose the $4$-angulation.}
\label{fig:angp}
\end{figure}
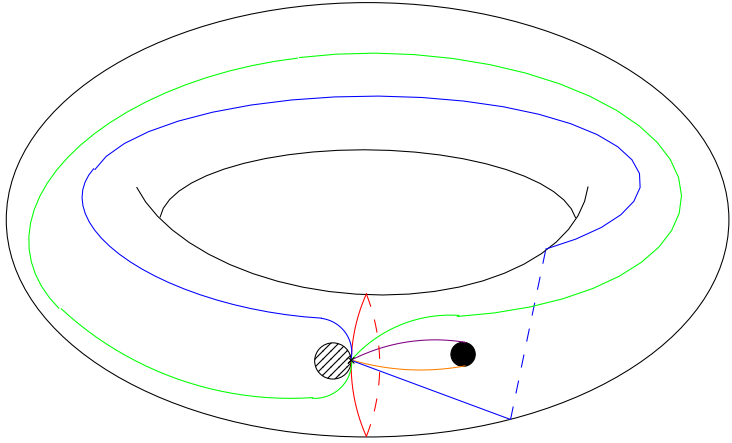

On figure \ref{fig:angp} there is another example of an $(m+2)$-angulation with one puncture. The drawing of the universal cover is left to the reader.

\begin{theo}\label{th:arcs}
Let $S$ be a surface. Let $g,b,c,p$ be respectively their gender, number of boundary components, number of marked points placed on the boundary components, and their number of punctures. We assume that:
\begin{enumerate}
\item Each boundary component contains at least one marked point, and the number of marked points is a multiple of $m$.
\item There is at least one boundary component.
\end{enumerate}

Then, the formula giving the number $n$ of arcs in an $(m+2)$-angulation of the surface is the following one :

\[ \frac{c}{m} + (1+\frac{2}{m})(b+2(g-1))+(2+\frac{1}{m})p=n. \]

Moreover, a necessary and sufficient condition to ensure the existence of an $(m+2)$-angulation is

\[c+2b+4g-4+p \equiv 0 [m]. \]

\end{theo}

\begin{proof}
It is immediate that the condition is necessary.
\vspace{10pt}

First step: We first consider the case where $p=0$:

Let us prove the result by induction on $(g,b,c)$ in the lexicographic order.

If $g=0$, $b=1$, then we are in the case where $S$ is associated with a quiver $Q$ of type $A_n$, and from the article of Baur and Marsh \cite{BM01} we know that there exists an $(m+2)$-angulation if and only if $c=(n+1)m+2$. Then we can check the result, and the condition is sufficient (indeed, $c+2-4=km$).

Let us now fix $(g,b,c)$ and assume that the result holds for all $(g',b',c') < (g,b,c)$. Let $S$ be a surface determined, up to homeomorphism by $(g,b,c)$. Let $n$ be the number of arcs necessary to an $(m+2)$-angulation. Let us consider two cases:

\vspace{10pt}

First case: $g \geq 1$:

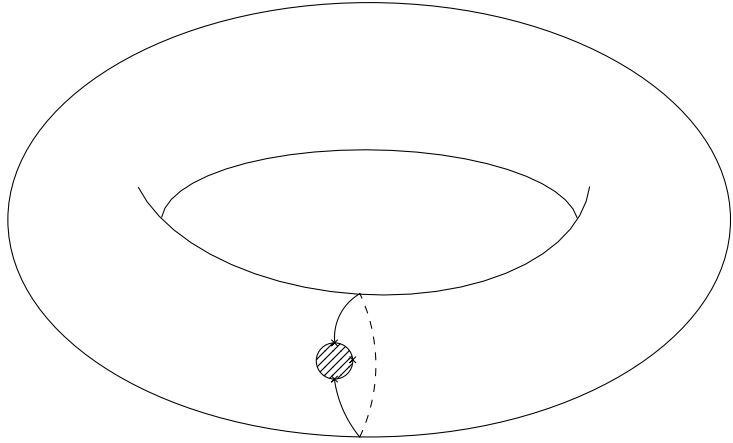
\begin{figure}[!h]
\centering
\begin{tikzpicture}[scale=0.8]
\draw [rotate around={0:(4.73,0)}] (4.73,0) ellipse (5.94cm and 3.6cm);
\draw [shift={(4.57,0.94)}] plot[domain=0.34:3.12,variable=\t]({-1*3.79*cos(\t r)+-0.09*2.16*sin(\t r)},{0.09*3.79*cos(\t r)+-1*2.16*sin(\t r)});
\draw [shift={(4.74,-0.07)}] plot[domain=0.11:3.08,variable=\t]({1*3.43*cos(\t r)+0.01*1.23*sin(\t r)},{-0.01*3.43*cos(\t r)+1*1.23*sin(\t r)});
\draw [fill=black,pattern=north east lines](4.16,-2.34) circle (0.3cm);
\draw [shift={(5,-1.95)}] plot[domain=2.09:3.25,variable=\t]({1*0.84*cos(\t r)+0*0.84*sin(\t r)},{0*0.84*cos(\t r)+1*0.84*sin(\t r)});
\draw [shift={(6.07,-2.38)}] plot[domain=3.28:3.83,variable=\t]({1*1.93*cos(\t r)+0*1.93*sin(\t r)},{0*1.93*cos(\t r)+1*1.93*sin(\t r)});
\draw [shift={(2.03,-2.41)},dash pattern=on 3pt off 3pt]  plot[domain=-0.44:0.44,variable=\t]({1*2.81*cos(\t r)+0*2.81*sin(\t r)},{0*2.81*cos(\t r)+1*2.81*sin(\t r)});
\begin{scriptsize}
\draw [color=black] (4.16,-2.04)-- ++(-1.5pt,-1.5pt) -- ++(3.0pt,3.0pt) ++(-3.0pt,0) -- ++(3.0pt,-3.0pt);
\draw [color=black] (4.16,-2.64)-- ++(-1.5pt,-1.5pt) -- ++(3.0pt,3.0pt) ++(-3.0pt,0) -- ++(3.0pt,-3.0pt);
\draw [color=black] (4.46,-2.32)-- ++(-1.5pt,-1.5pt) -- ++(3.0pt,3.0pt) ++(-3.0pt,0) -- ++(3.0pt,-3.0pt);
\end{scriptsize}
\end{tikzpicture}
\caption {Cutting along the arc decreases $g$ by $1$}
\label{fig:tore}
\end{figure}

There is at least one admissible arc going around the picture (as in figure \ref{fig:tore}). We will later see what an admissible arc is. If this is not the case, is means that $b \geq 2$ and we can go to the second case, and use an arc linking two boundary components (we know that there exists one).

Then we cut all along an admissible arc as in figure \ref{fig:tore}, decreasing $g$ by $1$ (There are $m$ marked points on the figure).
Then, if the new surface $S'$ is characterized by $(g',b',c')$, we have $g'=g-1$, $b'=b+1$, and $c'=c+2$.

Condition:
Let us suppose that 
$c+2b+4g-4=km, k \in {\mathbb{N}}$. Then, cutting along the same arc, we have the surface $S'$, characterized by $(g',b',c')$ and as $c+2+2(b+1)+4(g-1)-4=km$, we have also $c'+2b'+4g'-4=km$, and this is a sufficient condition to have an $(m+2)$-angulation of $S'$, given by $(g',b',c')$. Thus, the condition $c+2b+4g-4=km$ implies that we have an $(m+2)$-angulation of $S$ by adding an arc to the $(m+2)$-angulation of $S'$.

Any $(m+2)$-angulation of this surface $S'$ contains $n-1$ arcs. We also have that $(g,b',c') < (g,b,c)$ in the lexicographic order. By recurrence hypothesis we have
\[ \frac{c'}{m} + (1+\frac{2}{m})(b'+2g'-2)=n' \]

then,

\[ \frac{c}{m}+\frac{2}{m} + (1+\frac{2}{m})(b+1+2(g-1)-2)=n-1 \]

and so

\[ \frac{c}{m} + (1+\frac{2}{m})(b+2g-2)=n \]

and the result is proved.

\vspace{10pt}

Second case: $g=0$: cf figure \ref{fig:tore2}:

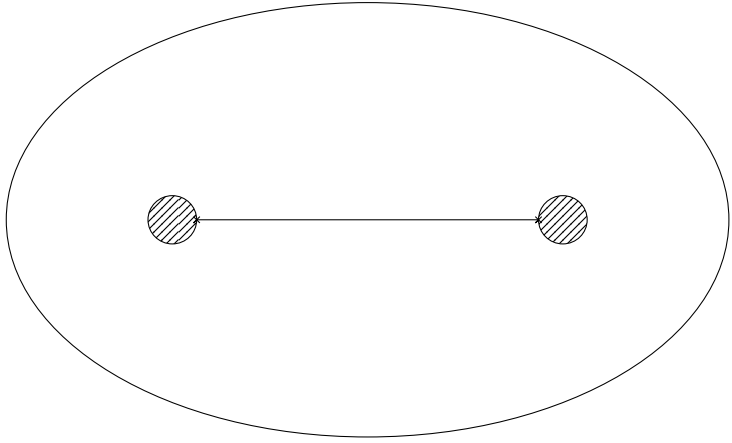
\begin{figure}[!h]
\centering
\begin{tikzpicture}[scale=0.8]
\draw [rotate around={0:(4.73,0)}] (4.73,0) ellipse (5.94cm and 3.6cm);
\draw [fill=black,pattern=north east lines](1.52,0) circle (0.4cm);
\draw [fill=black,pattern=north east lines](7.94,0) circle (0.4cm);
\draw (1.92,0)-- (7.54,0);
\begin{scriptsize}
\draw [color=black] (1.92,0)-- ++(-1.5pt,-1.5pt) -- ++(3.0pt,3.0pt) ++(-3.0pt,0) -- ++(3.0pt,-3.0pt);
\draw [color=black] (7.54,0)-- ++(-1.5pt,-1.5pt) -- ++(3.0pt,3.0pt) ++(-3.0pt,0) -- ++(3.0pt,-3.0pt);
\end{scriptsize}
\end{tikzpicture}
\caption{If $g=0$, we cut along an arc linking two boundary components}
\label{fig:tore2}
\end{figure}

Then, if $b=1$, the case has been treated in the initialisation. If $b \geq 2$, there is at least one marked point on each boundary component, thus we cut along an arc linking two boundary components.

Then, the new surface $S'$ is characterized by $(0,b',c')$ where $b'=b-1$, $c'=c+2$, and $n'=n-1$. We then have $(0,b',c') < (0,b,c)$. By recursion hypothesis, we have

\[ \frac{c'}{m} + (1+\frac{2}{m})(b'-2)=n' \]

then,

\[ \frac{c}{m} + \frac{2}{m} + (1 + \frac{2}{m})(b-3)=n-1 \]

and finally

\[ \frac{c}{m} + (1+\frac{2}{m})(b-2)=n \]

Condition:
Let us assume that $c+2b-4=km$ with $k \in {\mathbb{N}}$. Then, cutting along the same arc, we have the surface $S'$, characterized by $(0,b',c')$. Moreover,
$c+2+2(b-1)-4=km$ and $c'+2b'-4=km$ and this, by induction this condition is sufficient. Then, the condition $c+2b-4=km$ is sufficient to have an $(m+2)$-angulation of $S$, and this can be done by adding the cut arc to the surface $S'$ determined by $(0,b',c')$.

This ends the induction and the proof in the case where $p=0$.

\vspace{20pt}

Second step: $p$ is arbitrary:

Let us prove by induction on $p$ that the number $n$ of arcs for any $(m+2)$-angulation of the surface $S$, characterized by $(g,b,c,p)$ is given by the formula

\[ \frac{c}{m} +(1+\frac{2}{m} )(b+2g-2)+(2+\frac{1}{m})p=n \]

and a sufficient condition to be ensured of the existence of an $(m+2)$-angulation is $c+2b+4g-4+p=km$, where $k \in {\mathbb{N}}$.

The case where $p=0$ has been done at first step.

Let $S$ be a surface determined by $(g,b,c,p)$. Let us suppose the result for a number of punctures less than or equal to $p-1$. At the beginning, we have supposed $b \geq 1$. We also assume $p \geq 1$. Then there is at least one marked point on the boundary component. We cut along one arc linking a marked point to a puncture (cf figure \ref{fig:tore3}, the puncture is realized by a disk).

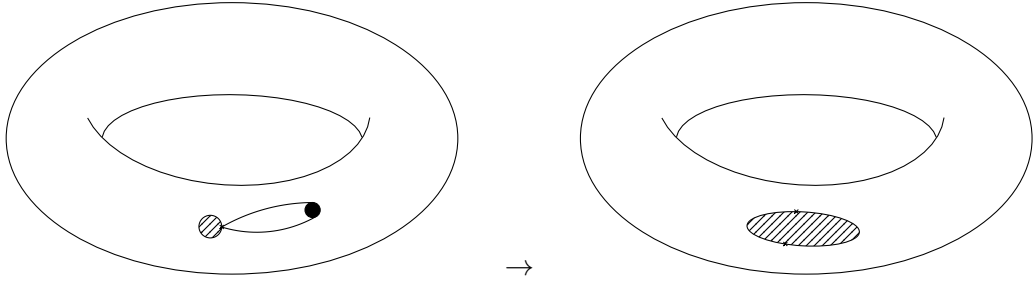
\begin{figure}[!h]
\centering
\begin{tikzpicture}[scale=0.5]
\draw [fill=black,fill opacity=1.0] (6.85,-1.9) circle (0.2cm);
\draw [rotate around={0:(4.73,0)}] (4.73,0) ellipse (5.94cm and 3.6cm);
\draw [shift={(4.57,0.94)}] plot[domain=0.34:3.12,variable=\t]({-1*3.79*cos(\t r)+-0.09*2.16*sin(\t r)},{0.09*3.79*cos(\t r)+-1*2.16*sin(\t r)});
\draw [shift={(4.74,-0.07)}] plot[domain=0.11:3.08,variable=\t]({1*3.43*cos(\t r)+0.01*1.23*sin(\t r)},{-0.01*3.43*cos(\t r)+1*1.23*sin(\t r)});
\draw [fill=black,pattern=north east lines](4.16,-2.34) circle (0.3cm);
\draw [shift={(6.89,-6.56)}] plot[domain=1.58:2.1,variable=\t]({1*4.86*cos(\t r)+0*4.86*sin(\t r)},{0*4.86*cos(\t r)+1*4.86*sin(\t r)});
\draw [shift={(5.38,0.64)}] plot[domain=4.41:5.23,variable=\t]({1*3.13*cos(\t r)+0*3.13*sin(\t r)},{0*3.13*cos(\t r)+1*3.13*sin(\t r)});
\begin{scriptsize}
\draw [color=black] (4.46,-2.35)-- ++(-1.5pt,-1.5pt) -- ++(3.0pt,3.0pt) ++(-3.0pt,0) -- ++(3.0pt,-3.0pt);
\end{scriptsize}
\end{tikzpicture}
\hspace{10pt}
$\to$
\hspace{10pt}
\begin{tikzpicture}[scale=0.5]
\draw [rotate around={0:(4.73,0)}] (4.73,0) ellipse (5.94cm and 3.6cm);
\draw [shift={(4.57,0.94)}] plot[domain=0.34:3.12,variable=\t]({-1*3.79*cos(\t r)+-0.09*2.16*sin(\t r)},{0.09*3.79*cos(\t r)+-1*2.16*sin(\t r)});
\draw [shift={(4.74,-0.07)}] plot[domain=0.11:3.08,variable=\t]({1*3.43*cos(\t r)+0.01*1.23*sin(\t r)},{-0.01*3.43*cos(\t r)+1*1.23*sin(\t r)});
\draw [rotate around={176.62:(4.65,-2.4)},fill=black,pattern=north east lines] (4.65,-2.4) ellipse (1.48cm and 0.45cm);
\begin{scriptsize}
\draw [color=black] (4.47,-1.94)-- ++(-1.5pt,-1.5pt) -- ++(3.0pt,3.0pt) ++(-3.0pt,0) -- ++(3.0pt,-3.0pt);
\draw [color=black] (4.18,-2.8)-- ++(-1.5pt,-1.5pt) -- ++(3.0pt,3.0pt) ++(-3.0pt,0) -- ++(3.0pt,-3.0pt);
\end{scriptsize}
\end{tikzpicture}

\caption{We cut the surface along the arcs $\alpha$}
\label{fig:tore3}
\end{figure}

Then we have a new surface $S'$, characterized by $(g',b',c',p')$ where $g'=g$, $b'=b$, $c'=c+1$, and $p'=p-1$. Note that $n'=n-2$.

By induction process,

\[ \frac{c'}{m} + (1+\frac{2}{m})(b'+2g'-2)+(2+\frac{1}{m})p'=n' \]

thus

\[ \frac{c}{m} + \frac{1}{m} + (1+ \frac{2}{m})(b+2g-2) + (2+ \frac{1}{m})(p-1)= n-2 \]

and so

\[ \frac{c}{m} + (1+\frac{2}{m})(b+2g'-2)+(2+\frac{1}{m})p=n \]

Condition:
Let us assume that $c+2b+4g-4+p=km$ where $k \in {\mathbb{N}}$. Then, cutting along the same arc, we have $S'$, characterized by $(g',b',c',p')$ and $c+1+2b+4g-4+p-1=km$, then $c'+2b'+4g'-4+p'=km$, and this, by induction is a sufficient condition to have an $(m+2)$-angulation on $S'$. Thus, the condition $c+2b+4g-4+p=km$ is sufficient.
This ends the induction and proves the theorem.
\end{proof}

\subsection{$(m+2)$-angulations and flips}

From here and all throughout the paper, weconsider that the marked surface has no punctures. Let us now watch a bit deeper into the $(m+2)$-angulations, and introduce the flip.

The following lemma is a result of Freedman, Hass, and Scott on general geodesics, which allows us to consider that the $m$-diagonals do not cross each other.

\begin{lem}\cite{FHS}
If $\Delta=\{\gamma_1,\cdots,\gamma_n \}$ is an $(m+2)$-angulation, then there exist $\{\alpha_1, \cdots, \alpha_n \}$ representatives of $\{\gamma_1,\cdots,\gamma_n \}$, such that for any $i$ and any $j \neq i$, $\alpha_i$ and $\alpha_j$ do not cross. This set is called a good set of representatives.
\end{lem}

We now introduce the twist which is an application defined on an arc of an $(m+2)$-angulation.

\begin{figure}[!h]
\centering
\begin{tikzpicture}[scale=0.8]
\draw (2.42,4.64)-- (3.84,3.26);
\draw (0.6,1.5)-- (2.3,-0.02);
\draw (2.42,4.64)-- (2.3,-0.02);
\draw (0.6,1.5)-- (3.84,3.26);
\begin{scriptsize}
\fill [color=black] (2.42,4.64) circle (1.5pt);
\draw[color=black] (2.58,4.9) node {$a$};
\fill [color=black] (2.3,-0.02) circle (1.5pt);
\draw[color=black] (2.46,0.24) node {$b$};
\draw[color=black] (2.5,3.84) node {$\alpha$};
\draw[color=black] (3.84,2.7) node {$\kappa\alpha$};
\end{scriptsize}
\end{tikzpicture}
\caption{Definition of a twist}
\label{fig:twist}
\end{figure}
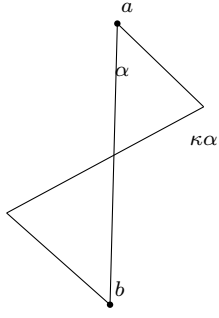

\begin{defi}
Let $\Delta$ be an $(m+2)$-angulation. Let $\alpha$ be an arc in it, linking the vertices $a$ and $b$. The twist of $\alpha$ in $\Delta$ is defined as follows:

Let $\alpha_a$ (respectively $\alpha_b$) be the side of the $(m+2)$-angle ending at $a$ (respectively at $b$) consecutive to $\alpha$ (respectively preceding $\alpha$). Then the twist of $\alpha$, namely $\kappa_{\Delta}(\alpha)$ is the path $\alpha_a \alpha \alpha_b$. See figure \ref{fig:twist} for an illustration of the twist.
\end{defi}

\begin{defi}
Consider $\Delta$ an $(m+2)$-angulation, it means that we choose $\{ \gamma_1,\cdots,\gamma_n \}$ good representatives. Let $\alpha$ be an arc in it. The flip of the $(m+2)$-angulation at $\alpha$ is defined by $\mu_\alpha\Delta=\Delta \setminus \{ \alpha \} \cup \{ \alpha^* \}$ where $\alpha^*$ is given by $\kappa_{\Delta}(\alpha)$, the twist of $\alpha$.
\end{defi}

\begin{rmk}
\begin{enumerate}
\item Note that the twist has an inverse, which consists in moving the arc counterclockwise. Then the flip is also invertible.
\item A flip does not change the number of arcs in the $(m+2)$-angulation.
\end{enumerate}
\end{rmk}

\section{Associating a graded quiver with superpotential with an $(m+2)$-angulation}

In this section, we associate a graded quiver and a superpotential with an $(m+2)$-angulation, and show in the unpunctured case the compatibility between the flip of an $(m+2)$-angulation and the mutation of the associated graded quiver with superpotential.

From now and all throughout the paper, let $n \geq 4$ and $m$ be positive integers. We set $S$ as a Riemann surface, in which there exist $(m+2)$-angulations.

\subsection{Definitions}
Let us now define the graded quiver with superpotential.

\begin{defi}\label{def:grquiv}
Let $\Delta$ be an $(m+2)$-angulation. Then the graded quiver with potential $(Q_\Delta^{gr},W_\Delta^{gr})$ is defined, up to quasi-isomorphism, as follows:
\begin{enumerate}
\item The vertices of $Q_\Delta^{gr}$ correspond to the arcs of the $(m+2)$-angulation. They are numbered from $1$ to $n$.
\item Between two vertices $i$ and $j$, there is an arrow when $i$ and $j$ form two sides of the same $(m+2)$-angle. If the number of sides from $i$ to $j$ is $r \in \{0,\ldots,m \}$, then the associated graduation is $r$.
\item The potential is given by the sum of all the $3$-cycles of degree $m-1$.
\end{enumerate}
\end{defi}

\begin{rmk}
Here, up to quasi-isomorphism means that we choose a representative quiver in the class of all quasi-isomorphic quivers with superpotential, using the previous definition.
\end{rmk}

\begin{prop}
There is an equivalent definition: the vertices are similarly defined, and for $i$ and $j$ two vertices, and $c$ an integer,
\[ 
q_{ij}^{r} = \left\{
    \begin{array}{ll}
        1 & \mbox{if } \kappa_\Delta^r(i) \mbox{ and } j \mbox{ share an oriented angle} \\
        0 & \mbox{otherwise.}
    \end{array}
    \right.
\] where $q_{ij}^{r}$ is the number of arrows from $i$ to $j$ of graduation $r$.
\end{prop}

\begin{proof}
We only have to show that the arrows are the same. If $i$ and $j$ form two sides of the polygon, with a graduation $r$, it means that if we apply the twist to $i$, then there will be $r-1$ edges between $ \kappa_\Delta(i)$ and $j$. Then if we apply the twist $r$ times, there will be no edge between $ \kappa_\Delta^r(i)$ and $j$, and they will share an oriented angle.

On the other hand, if $ \kappa_\Delta^r(i)$ and $j$ share an oriented angle, it suffices to apply the inverse of the twist $r$ times to make sure that $i$ and $j$ form two sides of a polygon, and that there are $r$ edges between $i$ and $j$.
\end{proof}

\begin{rmk}
Note that the quiver is "symmetric", which means that, if there is an arrow $\xymatrix{i \ar^r[r] & j}$, then there is an arrow $\xymatrix{j \ar^{m-r}[r] & i}$.
\end{rmk}

Let us draw an example. On figure \ref{fig:angu} there are three arcs forming a $4$-angulation, so the corresponding quiver has three vertices.

\begin{figure}[!h]
\centering
\begin{tikzpicture}[scale=0.6]
\draw [rotate around={0:(4.73,0)}] (4.73,0) ellipse (5.94cm and 3.6cm);
\draw [shift={(4.57,0.94)}] plot[domain=0.34:3.12,variable=\t]({-1*3.79*cos(\t r)+-0.09*2.16*sin(\t r)},{0.09*3.79*cos(\t r)+-1*2.16*sin(\t r)});
\draw [shift={(4.74,-0.07)}] plot[domain=0.11:3.08,variable=\t]({1*3.43*cos(\t r)+0.01*1.23*sin(\t r)},{-0.01*3.43*cos(\t r)+1*1.23*sin(\t r)});
\draw [fill=black,pattern=north east lines] (4.16,-2.34) circle (0.3cm);
\draw [shift={(5,-1.95)},color=red]  plot[domain=2.09:3.25,variable=\t]({1*0.84*cos(\t r)+0*0.84*sin(\t r)},{0*0.84*cos(\t r)+1*0.84*sin(\t r)});
\draw [shift={(6.07,-2.38)},color=red]  plot[domain=3.28:3.83,variable=\t]({1*1.93*cos(\t r)+0*1.93*sin(\t r)},{0*1.93*cos(\t r)+1*1.93*sin(\t r)});
\draw [shift={(2.03,-2.41)},dash pattern=on 3pt off 3pt,color=red]  plot[domain=-0.44:0.44,variable=\t]({1*2.81*cos(\t r)+0*2.81*sin(\t r)},{0*2.81*cos(\t r)+1*2.81*sin(\t r)});
\draw [color=green] (4.16,-2.04)-- (6.29,-1.09);
\draw [dash pattern=on 5pt off 5pt,color=green] (6.29,-1.09)-- (8.43,-2.82);
\draw [shift={(4.19,-0.23)},color=blue]  plot[domain=1.53:5.17,variable=\t]({-1*4.94*cos(\t r)+0.07*2.4*sin(\t r)},{-0.07*4.94*cos(\t r)+-1*2.4*sin(\t r)});
\draw [shift={(4.51,0.04)},color=blue]  plot[domain=-0.91:1.51,variable=\t]({-1*4.4*cos(\t r)+-0.03*2.09*sin(\t r)},{0.03*4.4*cos(\t r)+-1*2.09*sin(\t r)});
\draw [shift={(4.5,-0.48)},color=green]  plot[domain=2.32:5.6,variable=\t]({-0.99*5.38*cos(\t r)+0.11*3.74*sin(\t r)},{-0.11*5.38*cos(\t r)+-0.99*3.74*sin(\t r)});
\draw [shift={(4.47,0.78)},color=green]  plot[domain=-0.66:1.38,variable=\t]({-0.97*5.2*cos(\t r)+0.24*2.71*sin(\t r)},{-0.24*5.2*cos(\t r)+-0.97*2.71*sin(\t r)});
\begin{scriptsize}
\draw [color=black] (4.16,-2.04)-- ++(-1.5pt,-1.5pt) -- ++(3.0pt,3.0pt) ++(-3.0pt,0) -- ++(3.0pt,-3.0pt);
\draw [color=black] (4.16,-2.64)-- ++(-1.5pt,-1.5pt) -- ++(3.0pt,3.0pt) ++(-3.0pt,0) -- ++(3.0pt,-3.0pt);
\end{scriptsize}
\end{tikzpicture}
\caption{}
\label{fig:angu}
\end{figure}
We now draw the universal cover:
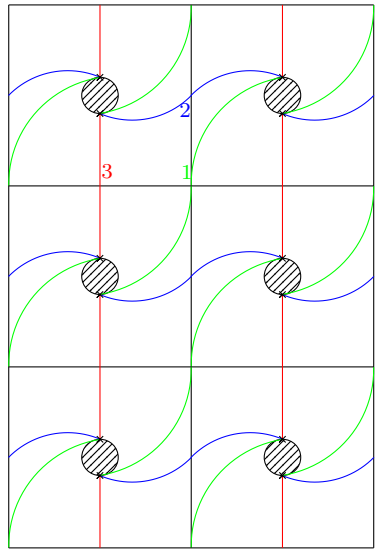
\begin{center}
\begin{tikzpicture}[scale=0.8]
\draw (-3,4)-- (-3,-5);
\draw (0,4)-- (0,-5);
\draw (3,-5)-- (3,4);
\draw (-3,4)-- (3,4);
\draw (-3,-5)-- (3,-5);
\draw (-3,1)-- (3,1);
\draw (-3,-2)-- (3,-2);
\draw [fill=black,pattern=north east lines] (-1.5,-0.5) circle (0.3cm);
\draw [fill=black,pattern=north east lines] (-1.5,-3.5) circle (0.3cm);
\draw [fill=black,pattern=north east lines] (-1.5,2.5) circle (0.3cm);
\draw [fill=black,pattern=north east lines] (1.5,2.5) circle (0.3cm);
\draw [fill=black,pattern=north east lines] (1.5,-0.5) circle (0.3cm);
\draw [fill=black,pattern=north east lines] (1.5,-3.5) circle (0.3cm);
\draw [shift={(-0.97,0.45)},color=blue]  plot[domain=4.31:5.51,variable=\t]({1*1.36*cos(\t r)+0*1.36*sin(\t r)},{0*1.36*cos(\t r)+1*1.36*sin(\t r)});
\draw [shift={(0.97,-1.45)},color=blue]  plot[domain=1.17:2.37,variable=\t]({1*1.36*cos(\t r)+0*1.36*sin(\t r)},{0*1.36*cos(\t r)+1*1.36*sin(\t r)});
\draw [shift={(-0.97,3.45)},color=blue]  plot[domain=4.31:5.51,variable=\t]({1*1.36*cos(\t r)+0*1.36*sin(\t r)},{0*1.36*cos(\t r)+1*1.36*sin(\t r)});
\draw [shift={(-0.97,-2.55)},color=blue]  plot[domain=4.31:5.51,variable=\t]({1*1.36*cos(\t r)+0*1.36*sin(\t r)},{0*1.36*cos(\t r)+1*1.36*sin(\t r)});
\draw [shift={(0.97,1.55)},color=blue]  plot[domain=1.17:2.37,variable=\t]({1*1.36*cos(\t r)+0*1.36*sin(\t r)},{0*1.36*cos(\t r)+1*1.36*sin(\t r)});
\draw [shift={(0.97,-4.45)},color=blue]  plot[domain=1.17:2.37,variable=\t]({1*1.36*cos(\t r)+0*1.36*sin(\t r)},{0*1.36*cos(\t r)+1*1.36*sin(\t r)});
\draw [shift={(2.03,3.45)},color=blue]  plot[domain=4.31:5.51,variable=\t]({1*1.36*cos(\t r)+0*1.36*sin(\t r)},{0*1.36*cos(\t r)+1*1.36*sin(\t r)});
\draw [shift={(2.03,0.45)},color=blue]  plot[domain=4.31:5.51,variable=\t]({1*1.36*cos(\t r)+0*1.36*sin(\t r)},{0*1.36*cos(\t r)+1*1.36*sin(\t r)});
\draw [shift={(2.03,-2.55)},color=blue]  plot[domain=4.31:5.51,variable=\t]({1*1.36*cos(\t r)+0*1.36*sin(\t r)},{0*1.36*cos(\t r)+1*1.36*sin(\t r)});
\draw [shift={(-2.03,-4.45)},color=blue]  plot[domain=1.17:2.37,variable=\t]({1*1.36*cos(\t r)+0*1.36*sin(\t r)},{0*1.36*cos(\t r)+1*1.36*sin(\t r)});
\draw [shift={(-2.03,-1.45)},color=blue]  plot[domain=1.17:2.37,variable=\t]({1*1.36*cos(\t r)+0*1.36*sin(\t r)},{0*1.36*cos(\t r)+1*1.36*sin(\t r)});
\draw [shift={(-2.03,1.55)},color=blue]  plot[domain=1.17:2.37,variable=\t]({1*1.36*cos(\t r)+0*1.36*sin(\t r)},{0*1.36*cos(\t r)+1*1.36*sin(\t r)});
\draw [color=red] (-1.5,2.8)-- (-1.5,4);
\draw [color=red] (-1.5,2.2)-- (-1.5,-0.2);
\draw [color=red] (-1.5,-0.8)-- (-1.5,-3.2);
\draw [color=red] (-1.5,-3.8)-- (-1.5,-5);
\draw [color=red] (1.5,2.8)-- (1.5,4);
\draw [color=red] (1.5,2.2)-- (1.5,-0.2);
\draw [color=red] (1.5,-0.8)-- (1.5,-3.2);
\draw [color=red] (1.5,-3.8)-- (1.5,-5);
\draw [shift={(-1.78,0.96)},color=green]  plot[domain=-1.41:0.02,variable=\t]({1*1.78*cos(\t r)+0*1.78*sin(\t r)},{0*1.78*cos(\t r)+1*1.78*sin(\t r)});
\draw [shift={(1.78,1.04)},color=green]  plot[domain=1.73:3.16,variable=\t]({1*1.78*cos(\t r)+0*1.78*sin(\t r)},{0*1.78*cos(\t r)+1*1.78*sin(\t r)});
\draw [shift={(-1.78,3.96)},color=green]  plot[domain=-1.41:0.02,variable=\t]({1*1.78*cos(\t r)+0*1.78*sin(\t r)},{0*1.78*cos(\t r)+1*1.78*sin(\t r)});
\draw [shift={(-1.78,-2.04)},color=green]  plot[domain=-1.41:0.02,variable=\t]({1*1.78*cos(\t r)+0*1.78*sin(\t r)},{0*1.78*cos(\t r)+1*1.78*sin(\t r)});
\draw [shift={(1.78,-1.96)},color=green]  plot[domain=1.73:3.16,variable=\t]({1*1.78*cos(\t r)+0*1.78*sin(\t r)},{0*1.78*cos(\t r)+1*1.78*sin(\t r)});
\draw [shift={(1.78,-4.96)},color=green]  plot[domain=1.73:3.16,variable=\t]({1*1.78*cos(\t r)+0*1.78*sin(\t r)},{0*1.78*cos(\t r)+1*1.78*sin(\t r)});
\draw [shift={(-1.22,1.04)},color=green]  plot[domain=1.73:3.16,variable=\t]({1*1.78*cos(\t r)+0*1.78*sin(\t r)},{0*1.78*cos(\t r)+1*1.78*sin(\t r)});
\draw [shift={(1.22,3.96)},color=green]  plot[domain=-1.41:0.02,variable=\t]({1*1.78*cos(\t r)+0*1.78*sin(\t r)},{0*1.78*cos(\t r)+1*1.78*sin(\t r)});
\draw [shift={(1.22,0.96)},color=green]  plot[domain=-1.41:0.02,variable=\t]({1*1.78*cos(\t r)+0*1.78*sin(\t r)},{0*1.78*cos(\t r)+1*1.78*sin(\t r)});
\draw [shift={(-1.22,-1.96)},color=green]  plot[domain=1.73:3.16,variable=\t]({1*1.78*cos(\t r)+0*1.78*sin(\t r)},{0*1.78*cos(\t r)+1*1.78*sin(\t r)});
\draw [shift={(1.22,-2.04)},color=green]  plot[domain=-1.41:0.02,variable=\t]({1*1.78*cos(\t r)+0*1.78*sin(\t r)},{0*1.78*cos(\t r)+1*1.78*sin(\t r)});
\draw [shift={(-1.22,-4.96)},color=green]  plot[domain=1.73:3.16,variable=\t]({1*1.78*cos(\t r)+0*1.78*sin(\t r)},{0*1.78*cos(\t r)+1*1.78*sin(\t r)});
\begin{scriptsize}
\draw [color=black] (-1.5,2.8)-- ++(-1.5pt,-1.5pt) -- ++(3.0pt,3.0pt) ++(-3.0pt,0) -- ++(3.0pt,-3.0pt);
\draw [color=black] (-1.5,2.2)-- ++(-1.5pt,-1.5pt) -- ++(3.0pt,3.0pt) ++(-3.0pt,0) -- ++(3.0pt,-3.0pt);
\draw [color=black] (-1.5,-0.2)-- ++(-1.5pt,-1.5pt) -- ++(3.0pt,3.0pt) ++(-3.0pt,0) -- ++(3.0pt,-3.0pt);
\draw [color=black] (-1.5,-0.8)-- ++(-1.5pt,-1.5pt) -- ++(3.0pt,3.0pt) ++(-3.0pt,0) -- ++(3.0pt,-3.0pt);
\draw [color=black] (-1.5,-3.2)-- ++(-1.5pt,-1.5pt) -- ++(3.0pt,3.0pt) ++(-3.0pt,0) -- ++(3.0pt,-3.0pt);
\draw [color=black] (-1.5,-3.8)-- ++(-1.5pt,-1.5pt) -- ++(3.0pt,3.0pt) ++(-3.0pt,0) -- ++(3.0pt,-3.0pt);
\draw [color=black] (1.5,2.8)-- ++(-1.5pt,-1.5pt) -- ++(3.0pt,3.0pt) ++(-3.0pt,0) -- ++(3.0pt,-3.0pt);
\draw [color=black] (1.5,2.2)-- ++(-1.5pt,-1.5pt) -- ++(3.0pt,3.0pt) ++(-3.0pt,0) -- ++(3.0pt,-3.0pt);
\draw [color=black] (1.5,-0.8)-- ++(-1.5pt,-1.5pt) -- ++(3.0pt,3.0pt) ++(-3.0pt,0) -- ++(3.0pt,-3.0pt);
\draw [color=black] (1.5,-3.8)-- ++(-1.5pt,-1.5pt) -- ++(3.0pt,3.0pt) ++(-3.0pt,0) -- ++(3.0pt,-3.0pt);
\draw [color=black] (1.5,-0.2)-- ++(-1.5pt,-1.5pt) -- ++(3.0pt,3.0pt) ++(-3.0pt,0) -- ++(3.0pt,-3.0pt);
\draw [color=black] (1.5,-3.2)-- ++(-1.5pt,-1.5pt) -- ++(3.0pt,3.0pt) ++(-3.0pt,0) -- ++(3.0pt,-3.0pt);
\draw[color=blue] (-0.1,2.26) node {$2$};
\draw[color=red] (-1.38,1.24) node {$3$};
\draw[color=green] (-0.07,1.22) node {$1$};
\end{scriptsize}
\end{tikzpicture}
\end{center}

and the corresponding graded quiver is given by

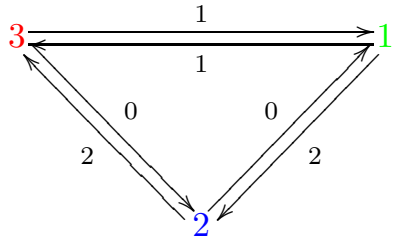
\[ \scalebox{1.2}{
\xymatrix@1 {
\textcolor{red}{3} \ar@<2pt>^1[rrrr] \ar@<2pt>^0[ddrr] & & & & \textcolor{green}{1} \ar@<2pt>^1[llll] \ar@<2pt>^2[ddll] \\
\\
& & \textcolor{blue}{2} \ar@<2pt>^2[uull] \ar@<2pt>^0[uurr] & & \\
}} \]

Let us now flip the $4$-angulation at vertex $2$ for instance. The new universal cover is the following one:

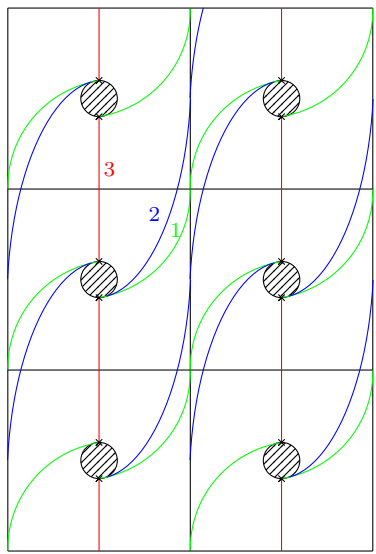
\begin{center}
\begin{tikzpicture}[scale=0.8]
\draw (-3,4)-- (-3,-5);
\draw (0,4)-- (0,-5);
\draw (3,-5)-- (3,4);
\draw (-3,4)-- (3,4);
\draw (-3,-5)-- (3,-5);
\draw (-3,1)-- (3,1);
\draw (-3,-2)-- (3,-2);
\draw [fill=black,pattern=north east lines] (-1.5,-0.5) circle (0.3cm);
\draw [fill=black,pattern=north east lines] (-1.5,-3.5) circle (0.3cm);
\draw [fill=black,pattern=north east lines] (-1.5,2.5) circle (0.3cm);
\draw [fill=black,pattern=north east lines] (1.5,2.5) circle (0.3cm);
\draw [fill=black,pattern=north east lines] (1.5,-0.5) circle (0.3cm);
\draw [fill=black,pattern=north east lines] (1.5,-3.5) circle (0.3cm);
\draw [color=red] (-1.5,2.8)-- (-1.5,4);
\draw [color=red] (-1.5,2.2)-- (-1.5,-0.2);
\draw [color=red] (-1.5,-0.8)-- (-1.5,-3.2);
\draw [color=red] (-1.5,-3.8)-- (-1.5,-5);
\draw [color=red] (1.5,2.8)-- (1.5,4);
\draw [color=red] (1.5,2.2)-- (1.5,-0.2);
\draw [color=red] (1.5,-0.8)-- (1.5,-3.2);
\draw [color=red] (1.5,-3.8)-- (1.5,-5);
\draw [shift={(-1.78,0.96)},color=green]  plot[domain=-1.41:0.02,variable=\t]({1*1.78*cos(\t r)+0*1.78*sin(\t r)},{0*1.78*cos(\t r)+1*1.78*sin(\t r)});
\draw [shift={(1.78,1.04)},color=green]  plot[domain=1.73:3.16,variable=\t]({1*1.78*cos(\t r)+0*1.78*sin(\t r)},{0*1.78*cos(\t r)+1*1.78*sin(\t r)});
\draw [shift={(-1.78,3.96)},color=green]  plot[domain=-1.41:0.02,variable=\t]({1*1.78*cos(\t r)+0*1.78*sin(\t r)},{0*1.78*cos(\t r)+1*1.78*sin(\t r)});
\draw [shift={(-1.78,-2.04)},color=green]  plot[domain=-1.41:0.02,variable=\t]({1*1.78*cos(\t r)+0*1.78*sin(\t r)},{0*1.78*cos(\t r)+1*1.78*sin(\t r)});
\draw [shift={(1.78,-1.96)},color=green]  plot[domain=1.73:3.16,variable=\t]({1*1.78*cos(\t r)+0*1.78*sin(\t r)},{0*1.78*cos(\t r)+1*1.78*sin(\t r)});
\draw [shift={(1.78,-4.96)},color=green]  plot[domain=1.73:3.16,variable=\t]({1*1.78*cos(\t r)+0*1.78*sin(\t r)},{0*1.78*cos(\t r)+1*1.78*sin(\t r)});
\draw [shift={(-1.22,1.04)},color=green]  plot[domain=1.73:3.16,variable=\t]({1*1.78*cos(\t r)+0*1.78*sin(\t r)},{0*1.78*cos(\t r)+1*1.78*sin(\t r)});
\draw [shift={(1.22,3.96)},color=green]  plot[domain=-1.41:0.02,variable=\t]({1*1.78*cos(\t r)+0*1.78*sin(\t r)},{0*1.78*cos(\t r)+1*1.78*sin(\t r)});
\draw [shift={(1.22,0.96)},color=green]  plot[domain=-1.41:0.02,variable=\t]({1*1.78*cos(\t r)+0*1.78*sin(\t r)},{0*1.78*cos(\t r)+1*1.78*sin(\t r)});
\draw [shift={(-1.22,-1.96)},color=green]  plot[domain=1.73:3.16,variable=\t]({1*1.78*cos(\t r)+0*1.78*sin(\t r)},{0*1.78*cos(\t r)+1*1.78*sin(\t r)});
\draw [shift={(1.22,-2.04)},color=green]  plot[domain=-1.41:0.02,variable=\t]({1*1.78*cos(\t r)+0*1.78*sin(\t r)},{0*1.78*cos(\t r)+1*1.78*sin(\t r)});
\draw [shift={(-1.22,-4.96)},color=green]  plot[domain=1.73:3.16,variable=\t]({1*1.78*cos(\t r)+0*1.78*sin(\t r)},{0*1.78*cos(\t r)+1*1.78*sin(\t r)});
\draw [shift={(-1.46,-0.32)},color=blue]  plot[domain=0.13:1.53,variable=\t]({-0.03*3.48*cos(\t r)+1*1.47*sin(\t r)},{-1*3.48*cos(\t r)+-0.03*1.47*sin(\t r)});
\draw [shift={(1.46,-0.68)},color=blue]  plot[domain=0.13:1.53,variable=\t]({0.03*3.48*cos(\t r)+-1*1.47*sin(\t r)},{1*3.48*cos(\t r)+0.03*1.47*sin(\t r)});
\draw [shift={(-1.46,2.68)},color=blue]  plot[domain=0.13:1.53,variable=\t]({-0.03*3.48*cos(\t r)+1*1.47*sin(\t r)},{-1*3.48*cos(\t r)+-0.03*1.47*sin(\t r)});
\draw [shift={(1.46,-3.68)},color=blue]  plot[domain=0.13:1.53,variable=\t]({0.03*3.48*cos(\t r)+-1*1.47*sin(\t r)},{1*3.48*cos(\t r)+0.03*1.47*sin(\t r)});
\draw [shift={(1.54,2.68)},color=blue]  plot[domain=0.13:1.53,variable=\t]({-0.03*3.48*cos(\t r)+1*1.47*sin(\t r)},{-1*3.48*cos(\t r)+-0.03*1.47*sin(\t r)});
\draw [shift={(1.54,-0.32)},color=blue]  plot[domain=0.13:1.53,variable=\t]({-0.03*3.48*cos(\t r)+1*1.47*sin(\t r)},{-1*3.48*cos(\t r)+-0.03*1.47*sin(\t r)});
\draw [shift={(-1.54,-3.68)},color=blue]  plot[domain=0.13:1.53,variable=\t]({0.03*3.48*cos(\t r)+-1*1.47*sin(\t r)},{1*3.48*cos(\t r)+0.03*1.47*sin(\t r)});
\draw [shift={(-1.54,-0.68)},color=blue]  plot[domain=0.13:1.53,variable=\t]({0.03*3.48*cos(\t r)+-1*1.47*sin(\t r)},{1*3.48*cos(\t r)+0.03*1.47*sin(\t r)});
\draw [shift={(1.46,2.32)},color=blue]  plot[domain=4.22:4.67,variable=\t]({-0.03*3.48*cos(\t r)+1*1.47*sin(\t r)},{-1*3.48*cos(\t r)+-0.03*1.47*sin(\t r)});
\begin{scriptsize}
\draw [color=black] (-1.5,2.8)-- ++(-1.5pt,-1.5pt) -- ++(3.0pt,3.0pt) ++(-3.0pt,0) -- ++(3.0pt,-3.0pt);
\draw [color=black] (-1.5,2.2)-- ++(-1.5pt,-1.5pt) -- ++(3.0pt,3.0pt) ++(-3.0pt,0) -- ++(3.0pt,-3.0pt);
\draw [color=black] (-1.5,-0.2)-- ++(-1.5pt,-1.5pt) -- ++(3.0pt,3.0pt) ++(-3.0pt,0) -- ++(3.0pt,-3.0pt);
\draw [color=black] (-1.5,-0.8)-- ++(-1.5pt,-1.5pt) -- ++(3.0pt,3.0pt) ++(-3.0pt,0) -- ++(3.0pt,-3.0pt);
\draw [color=black] (-1.5,-3.2)-- ++(-1.5pt,-1.5pt) -- ++(3.0pt,3.0pt) ++(-3.0pt,0) -- ++(3.0pt,-3.0pt);
\draw [color=black] (-1.5,-3.8)-- ++(-1.5pt,-1.5pt) -- ++(3.0pt,3.0pt) ++(-3.0pt,0) -- ++(3.0pt,-3.0pt);
\draw [color=black] (1.5,2.8)-- ++(-1.5pt,-1.5pt) -- ++(3.0pt,3.0pt) ++(-3.0pt,0) -- ++(3.0pt,-3.0pt);
\draw [color=black] (1.5,2.2)-- ++(-1.5pt,-1.5pt) -- ++(3.0pt,3.0pt) ++(-3.0pt,0) -- ++(3.0pt,-3.0pt);
\draw [color=black] (1.5,-0.8)-- ++(-1.5pt,-1.5pt) -- ++(3.0pt,3.0pt) ++(-3.0pt,0) -- ++(3.0pt,-3.0pt);
\draw [color=black] (1.5,-3.8)-- ++(-1.5pt,-1.5pt) -- ++(3.0pt,3.0pt) ++(-3.0pt,0) -- ++(3.0pt,-3.0pt);
\draw[color=red] (-1.33,1.33) node {$3$};
\draw [color=black] (1.5,-0.2)-- ++(-1.5pt,-1.5pt) -- ++(3.0pt,3.0pt) ++(-3.0pt,0) -- ++(3.0pt,-3.0pt);
\draw [color=black] (1.5,-3.2)-- ++(-1.5pt,-1.5pt) -- ++(3.0pt,3.0pt) ++(-3.0pt,0) -- ++(3.0pt,-3.0pt);
\draw[color=green] (-0.24,0.32) node {$1$};
\draw[color=blue] (-0.59,0.59) node {$2$};
\end{scriptsize}
\end{tikzpicture}
\end{center}

Now we build from this new $(m+2)$-angulation the associated quiver:

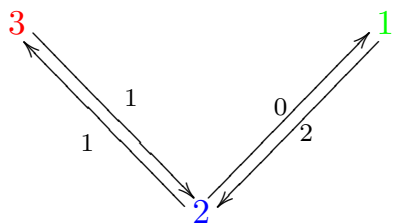
\[ \scalebox{1.2}{
\xymatrix@1 {
\textcolor{red}{3}  \ar@<2pt>^1[ddrr] & & & & \textcolor{green}{1} \ar@<2pt>_0[ddll] \\
\\
& & \textcolor{blue}{2} \ar@<2pt>^1[uull] \ar@<2pt>_2[uurr] & & \\
}} \]

The next section states that this new graded quiver with potential is quasi-isomorphic to the graded mutation in the sense of Oppermann of the former graded quiver with superpotential. In this example, the mutation of the first quiver gives this one:

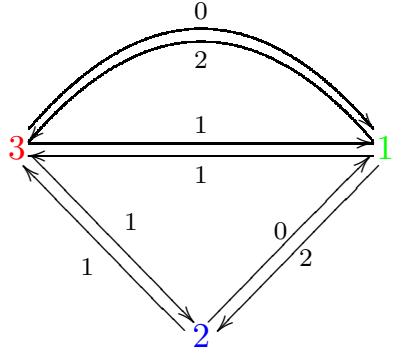
\[ \scalebox{1.2}{
\xymatrix@1 {
\textcolor{red}{3} \ar@<2pt>^1[rrrr]  \ar@<2pt>@/^3pc/^0[rrrr] \ar@<2pt>^1[ddrr] & & & & \textcolor{green}{1} \ar@<2pt>^1[llll] \ar@<2pt>@/_3pc/^2[llll] \ar@<2pt>_0[ddll] \\
\\
& & \textcolor{blue}{2} \ar@<2pt>^1[uull] \ar@<2pt>_2[uurr] & & \\
}} \]

which turns to be quasi-isomorphic to the quiver above.

\begin{prop}
Let $\Delta$ be an $(m+2)$-angulation of $S$. Let $(Q_\Delta,W_\Delta)$ be the associated quiver with superpotential. Then, this quiver satisfies the hypotheses of Guo and Oppermann, it means that:
\begin{itemize}
\item The graduations belong to $\{0,\cdots,m\}$.
\item If $\xymatrix{ i\ar^r[r] & j}$ and $\xymatrix{ j\ar^{r'}[r] & i}$, then $r'=m-r$.
\item $W$ has degree $m-1$.
\item The necklace bracket $\{W,W\}$ vanishes.
\end{itemize}
\end{prop}

\begin{proof}
\begin{itemize}
\item Let $i$ and $j$ be two vertices sharing an $(m+2)$-angle. Then, there can be at most $m$ sides between them ($i$ and $j$ forming both missing edges to build an $(m+2)$-angle.
\item If there are $r$ edges from $i$ to $j$ clockwise and $r'$ from $j$ to $i$, then $r+r'+2=m+2$ (we add $2$ because $i$ and $j$ belong to the $(m+2)$-angle). Then $r'=m-r$.
\item By definition, the degree of a potential is the sum of the degrees of the cycles forming it.
\item The superpotential $W$ only depends on the arrows. Then its necklace bracket vanishes, as we can see from Van Den Bergh in \cite[page 2]{VdB}
\end{itemize}
\end{proof}

\subsection{Compatibility with mutation in the unpunctured case}

\vspace{20pt}

Let us now draw another issue: We take an $(m+2)$-angulation of a surface, associate the corresponding graded quiver with superpotential, mutate it at a vertex $i$, flip the $(m+2)$-angulation at $i$, and analyse what happens.

\begin{theo}\label{th:comp}
Let $(S,M)$ be an unpunctured surface admitting $(m+2)$-angulations. Let $\Delta$ be an $(m+2)$-angulation. Let $(Q_\Delta,W_\Delta)$ be the associated graded quiver with potential. Let $i$ be an arc in $\Delta$. We call by $\Delta'=\mu_i(\Delta)$ the new $(m+2)$-angulation obtained from $\Delta$ by flipping at $i$.

We also call by $\mu_i(Q_\Delta,W_\Delta)$ the graded mutation in the sense of Oppermann of the graded quiver with superpotential $(Q_\Delta,W_\Delta)$. Then we have:

\[ \mu_i(Q_\Delta,W_\Delta)=(Q_{\mu_i(\Delta)},W_{\mu_i(\Delta)}). \]
\end{theo}

\begin{proof}
Note that the vertices are the same. We have to check the equality with respect to the arrows, graduations and superpotential.

Let $j$ and $k$ be two arcs in $\Delta$ and $\Delta'=\mu_i(\Delta)$. We denote by $q_{jk}^r$ the number of arrows frolm $j$ to $k$ of graduation $r$ in $Q_\Delta$. We denote by $\overline{q}_{jk}^r$ the number of arrows frolm $j$ to $k$ of graduation $r$ in $Q_{\mu_i(\Delta}$. We denote by $\tilde{q}_{jk}^r$ the number of arrows frolm $j$ to $k$ of graduation $r$ in $\mu_i(Q_\Delta)$.

The aim is to show in all cases, that $\tilde{q}_{jk}^r=\overline{q}_{jk}^r$.

\begin{itemize}
\item If $k=i$.
\begin{itemize}
\item If $r \neq m$.
Then $\tilde{q}_{ji}^r=q_{ji}^{r-1}$ by definition of the quiver mutation. Moreover, the flip of the arc put the difference from $j$ to $i$ from $r$ to $r-1$. Then $\overline{q}_{ji}^r=q_{ji}^{r-1}$. See figure \ref{fig:mutinm} for an explanation.

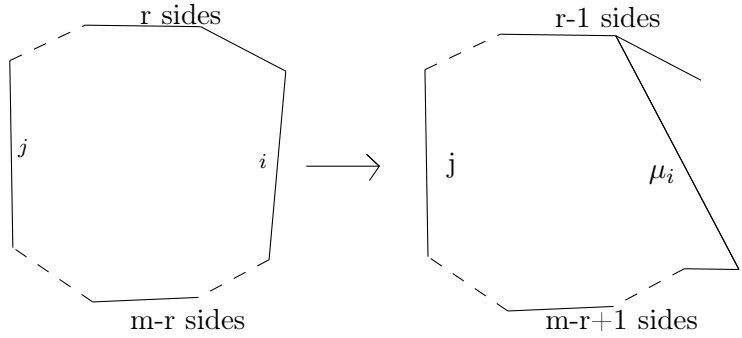
\begin{figure}[!h]
\centering
\begin{tikzpicture}[scale=0.5]
\draw [dash pattern=on 5pt off 5pt] (-1.24,1.92)-- (0.74,2.86);
\draw (0.74,2.86)-- (3.78,2.82);
\draw (3.78,2.82)-- (6.02,1.64);
\draw (6.02,1.64)-- (5.58,-3.36);
\draw [dash pattern=on 5pt off 5pt] (5.58,-3.36)-- (3.72,-4.36);
\draw (3.72,-4.36)-- (0.94,-4.48);
\draw [dash pattern=on 5pt off 5pt] (0.94,-4.48)-- (-1.16,-3.04);
\draw (-1.16,-3.04)-- (-1.24,1.92);
\draw (1.92,3.7) node[anchor=north west] {r sides};
\draw (1.66,-4.4) node[anchor=north west] {m-r sides};
\draw [dash pattern=on 5pt off 5pt] (9.68,1.68)-- (11.66,2.62);
\draw (11.66,2.62)-- (14.7,2.58);
\draw (14.7,2.58)-- (16.94,1.4);
\draw [dash pattern=on 5pt off 5pt] (16.5,-3.6)-- (14.64,-4.6);
\draw (14.64,-4.6)-- (11.86,-4.72);
\draw [dash pattern=on 5pt off 5pt] (11.86,-4.72)-- (9.76,-3.28);
\draw (9.76,-3.28)-- (9.68,1.68);
\draw (12.84,3.6) node[anchor=north west] {r-1 sides};
\draw (12.58,-4.4) node[anchor=north west] {m-r+1 sides};
\draw (10,-0.24) node[anchor=north west] {j};
\draw (15.3,-0.6) node[anchor=north west] {$\mu_i$};
\draw (14.7,2.58)-- (17.94,-3.62);
\draw (16.5,-3.6)-- (17.94,-3.62);
\draw (17.94,-3.62)-- (14.7,2.58);
\draw (6.56,-0.88)-- (8.48,-0.88);
\draw (8.48,-0.88)-- (8.14,-0.56);
\draw (8.48,-0.88)-- (8.22,-1.26);
\draw (8.14,-0.56)-- (8.48,-0.88);
\begin{scriptsize}
\draw[color=black] (5.42,-0.74) node {$i$};
\draw[color=black] (-0.88,-0.4) node {$j$};
\end{scriptsize}
\end{tikzpicture}
\caption{Mutation in the case $k=i$ and $m \neq r$}
\label{fig:mutinm}
\end{figure}
\item If $r=m$. Then we have an arrow from $j$ to $i$ of graduation $m$, which becomes in the mutated quiver an arrow from $j$ to $i$ of graduation $0$ (and an inverse arrow from $i$ to $j$ of graduation $m$). Then $\tilde{q}_{ji}^0=q_{ji}^{m}$. Moreover, we are in the situation of figure \ref{fig:mutim}

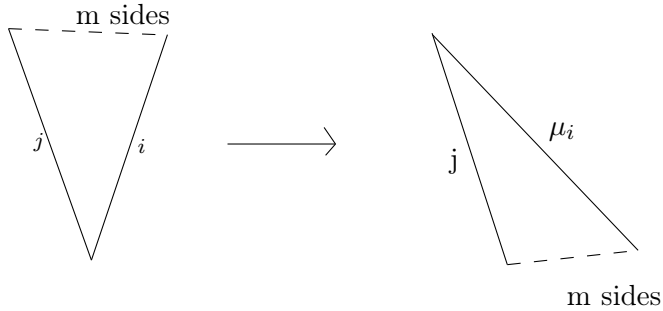
\begin{figure}[!h]
\centering
\begin{tikzpicture}[scale=0.5]
\draw (-0.2,2.32)-- (1.98,-3.82);
\draw (1.98,-3.82)-- (3.98,2.16);
\draw (10.94,2.2)-- (12.92,-3.94);
\draw (5.56,-0.74)-- (8.4,-0.74);
\draw (8.4,-0.74)-- (8.1,-0.3);
\draw (8.4,-0.74)-- (8.14,-1.04);
\draw (10.96,2.15)-- (16.36,-3.56);
\draw [dash pattern=on 5pt off 5pt] (-0.2,2.32)-- (3.98,2.15);
\draw [dash pattern=on 5pt off 5pt] (12.92,-3.94)-- (16.36,-3.56);
\draw (1.3,3.24) node[anchor=north west] {m sides};
\draw (14.22,-4.2) node[anchor=north west] {m sides};
\draw (11.14,-0.6) node[anchor=north west] {j};
\draw (13.74,0.1) node[anchor=north west] {$\mu_i$};
\begin{scriptsize}
\draw[color=black] (0.6,-0.7) node {$j$};
\draw[color=black] (3.3,-0.78) node {$i$};
\end{scriptsize}
\end{tikzpicture}
\caption{Mutation in the case $k=i$ and $m = r$}
\label{fig:mutim}
\end{figure}
Then, $\overline{q}_{ji}^0=q_{ji}^{m}$ and this concludes the case $i=k$.
\end{itemize}
\item If $i=j$.
By definition of $Q_{\mu_i(\Delta)}$, we have $\overline{q}_{ik}^r=\overline{q}_{ki}^{m-r}$. In addition, by definition of the graded quiver mutation, we have $\tilde{q}_{ik}^r=\tilde{q}_{ki}^{m-r}$. Then we come back to the previous case to conclude.
\item Else, assume that, in $Q_\Delta$, we have $\xymatrix@1{i\ar[r]^{r} & k\ar[r]^{0} & j}$ and $i$ and $j$ are not two sides of the same $(m+2)$-angle. Then $\overline{q}_{ij}^{r}=q_{ij}^{r}+1$. We notice that we cannot have two arrows from $k$ to $j$ of graduation $0$. We are in a situation of the following type:

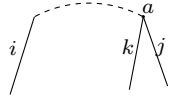
\begin{center}
\begin{tikzpicture}[scale=0.4]
\draw (1.26,2.8)-- (0.46,0.22);
\draw [shift={(3.05,-0.51)},dash pattern=on 2pt off 2pt]  plot[domain=1.08:2.07,variable=\t]({1*3.76*cos(\t r)+0*3.76*sin(\t r)},{0*3.76*cos(\t r)+1*3.76*sin(\t r)});
\draw (4.84,2.8)-- (4.38,0.38);
\draw (4.84,2.8)-- (5.64,0.5);
\begin{scriptsize}
\draw[color=black] (0.56,1.76) node {$i$};
\fill [color=black] (4.84,2.8) circle (1.5pt);
\draw[color=black] (5,3.06) node {$a$};
\draw[color=black] (4.34,1.8) node {$k$};
\draw[color=black] (5.42,1.78) node {$j$};
\end{scriptsize}
\end{tikzpicture}
\end{center}

Let then $a$ be the common vertex of $k$ and $j$. Let $b$ be the other vertex of $j$ (where $a$ or $b$ can be central polygons). Then $\kappa_\Delta(k)$ and $j$ have $b$ as a common vertex. As, in $\Delta$, $\kappa_\Delta^c(i)$ and $k$ have one common vertex $a$, then $\kappa_{\mathcal{T}}^{r}(i)$ and $j$ have one common vertex. As a consequence, $\tilde{q}_{ij}^{r}=q_{ij}^{r}+1$. There is no loss because there are not any arrows from $i$ to $j$ of graduation different from $c$. We note that there can be two arrows of different graduations appearing from $j$ to $k$, but as this is forbidden, we remove them, and this confirms that there are no two such arcs in $\mu_k(\Delta)$. This case is symmetric to the case where $i$ and $k$ share a common vertex.

Assume that we are in a situation of the following type:

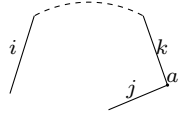
\begin{center}
\begin{tikzpicture}[scale=0.4]
\draw (1.26,2.8)-- (0.46,0.22);
\draw [shift={(3.05,-0.51)},dash pattern=on 2pt off 2pt]  plot[domain=1.08:2.07,variable=\t]({1*3.76*cos(\t r)+0*3.76*sin(\t r)},{0*3.76*cos(\t r)+1*3.76*sin(\t r)});
\draw (4.84,2.8)-- (5.64,0.5);
\draw (5.64,0.5)-- (3.7,-0.32);
\begin{scriptsize}
\draw[color=black] (0.56,1.76) node {$i$};
\draw[color=black] (5.46,1.78) node {$k$};
\fill [color=black] (5.64,0.5) circle (1.5pt);
\draw[color=black] (5.8,0.76) node {$a$};
\draw[color=black] (4.46,0.38) node {$j$};
\end{scriptsize}
\end{tikzpicture}
\end{center}

It means that we have a quiver:
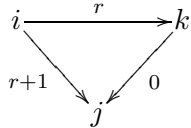
\[ \xymatrix@1{
i\ar[rr]^{r} \ar[dr]_{r+1} & & k \ar[dl]^{0} \\
& j
} \]

Then $\overline{q}_{ij}^{r}=q_{ij}^{r}-1=0$, because the arrows get erased.
Moreover, under the notations of the third case, $\kappa_\Delta(k)$ and $j$ have $b$ as a common vertex, and $i$ and $j$ do not share a common polygon anymore. Then $\tilde{q}_{ij}^{r}=q_{ij}^{r}-1$.
\end{itemize}

\begin{rmk}
The case of the Oppermann mutation of graded quivers is slightly different from the case of mutation of colored quivers introduced by Buan and Thomas (in \cite{BT}). Indeed, if two arrows $\alpha$, $\beta$ of different graduation appear from $i$ to $j$ (which is forbidden in Buan and Thomas' work), it suffices to factorize by $\alpha \beta^{-1}$.

For instance, see section 4 in Oppermann's work in \cite{O}, example 4.1. He explains that a quiver with two arrows of different graduation is quasi isomorphic to the same quiver where both arrows have been erased.
\end{rmk}

Now it remains to show that the potentials of both quivers are the same in order to finish the proof.

Let $W$ be the sum of all $3$-cycles of degree $m-1$ in $\Delta$. Let $\tilde{W}$ be the mutated potential, and let $\overline{W}$ be the sum of all $3$-cycles of degree $m-1$.

We show that $\tilde{W}=\overline{W}$.

\begin{itemize}
\item All the $3$-cycle in which $i$ does not lie are in $\tilde{W}$ (because $\tilde{W}=[W]+\sum \alpha \varphi \varphi^{\text{op}}\alpha^*$, and these cycles lie in $[W]$) and $\overline{W}$ (by definition).
\item Now, if there is a cycle in which $i$ appears. We are in this situation:
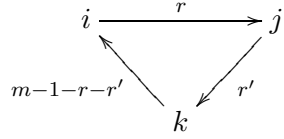
\[\xymatrix{i \ar^r[rr] & & j \ar^{r'}[dl] \\ & k \ar^{m-1-r-r'}[ul] &}\]

Let us decompose $W$ in $W=W'+W''$ where $W'$ contains all the cycles in which $i$ does not appear, and $W''$ contains all the cycles in which $i$ appears.
\begin{itemize}
\item If $r=0$. Then the new quiver $Q'$ is of this type:
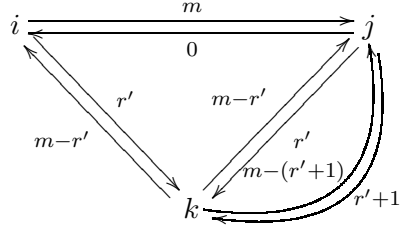
\[\xymatrix{i \ar@<0.5ex>^m[rrrr] \ar@<0.5ex>^{r'}[ddrr] & & & & j \ar@<0.5ex>^0[llll] \ar@<0.5ex>^{r'}[ddll] \ar@/^3pc/@<0.5ex>_{m-(r'+1)}[ddll] \\ & & & & \\ & & k \ar@<0.5ex>^{m-r'}[uull] \ar@<0.5ex>^{m-r'}[uurr] \ar@/_3pc/@<0.5ex>_{r'+1}[uurr] & &}\]

This quiver, having two arrows from $j$ to $k$ of different graduations, is quasi isomorphic to:
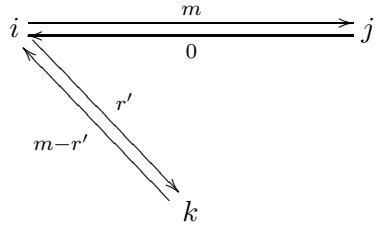
\[\xymatrix{i \ar@<0.5ex>^m[rrrr] \ar@<0.5ex>^{r'}[ddrr] & & & & j \ar@<0.5ex>^0[llll] \\ & & & & \\ & & k \ar@<0.5ex>^{m-r'}[uull] & &}\]

Then the cycle in which $i$ lies disappears. Then $\overline{W}=W'$.

Now we calculate $\tilde{W}$. By definition, \[\tilde{W}=[W]+\sum \alpha \varphi \varphi^{\text{op}}\alpha^*\] where the sum lies for $\alpha:i \to j$ of graduation $0$, and $\varphi:j \to i$ of graduation $r$ in $Q$.

The previous quiver was:
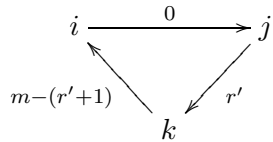
\[\xymatrix{i \ar^0[rr] & & j \ar^{r'}[dl] \\ & k \ar^{m-(r'+1)}[ul] &}\]
It becomes:
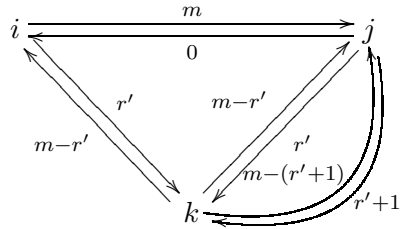
\[\xymatrix{i \ar@<0.5ex>^m[rrrr] \ar@<0.5ex>^{r'}[ddrr] & & & & j \ar@<0.5ex>^0[llll] \ar@<0.5ex>^{r'}[ddll] \ar@/^3pc/@<0.5ex>_{m-(r'+1)}[ddll] \\ & & & & \\ & & k \ar@<0.5ex>^{m-r'}[uull] \ar@<0.5ex>^{m-r'}[uurr] \ar@/_3pc/@<0.5ex>_{r'+1}[uurr] & &}\]

Then, $\tilde{W}=W'+\alpha \varphi \varphi^{\text{op}}\alpha^*$. But the second member of the sum contains a $2$-cycle, which is zero up to quasi isomorphism.
\item If $r \neq 0$. We have the quiver:
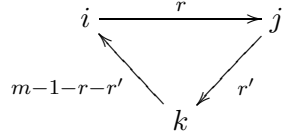
\[\xymatrix{i \ar^r[rr] & & j \ar^{r'}[dl] \\ & k \ar^{m-1-r-r'}[ul] &}\]
It becomes:
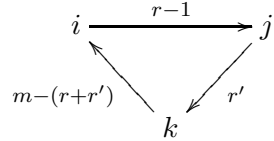
\[\xymatrix{i \ar^{r-1}[rr] & & j \ar^{r'}[dl] \\ & k \ar^{m-(r+r')}[ul] &}\]
This cycle (called $w$) is of degree $m-(r+r')+r-1+r'=m-1$. Then, $\overline{W}=W'+w$. In $W$, the cycle $w$ already appeared in $W''$ (since $i$ is involved). On the other hand \[\tilde{W}=[W]+\sum \alpha \varphi \varphi^{\text{op}}\alpha^*\] where the sum lies for $\alpha:i \to j$ of graduation $0$, and $\varphi:j \to i$ of graduation $r$ in $Q$. However, there is no arrow $i \to j$ of graduation $0$ since we have supposed $r \neq 0$.
Then $\overline{W}=\tilde{W}=W$.
\end{itemize}
\end{itemize}
\end{proof}

Let us draw an example:

Let $m=2$. We consider the following $(m+2)$-angulation in figure \ref{fig:ex}:

\begin{figure}[!h]
\centering
\begin{tikzpicture}[scale=0.5]
\fill[fill=black,fill opacity=0.1] (-0.12,1.88) -- (2.2,1.88) -- (4.08,3.24) -- (4.79,5.45) -- (4.08,7.66) -- (2.2,9.02) -- (-0.12,9.02) -- (-2,7.66) -- (-2.71,5.45) -- (-2,3.24) -- cycle;
\draw (-0.12,1.88)-- (2.2,1.88);
\draw (2.2,1.88)-- (4.08,3.24);
\draw (4.08,3.24)-- (4.79,5.45);
\draw (4.79,5.45)-- (4.08,7.66);
\draw (4.08,7.66)-- (2.2,9.02);
\draw (2.2,9.02)-- (-0.12,9.02);
\draw (-0.12,9.02)-- (-2,7.66);
\draw (-2,7.66)-- (-2.71,5.45);
\draw (-2.71,5.45)-- (-2,3.24);
\draw (-2,3.24)-- (-0.12,1.88);
\draw (-0.12,9.02)-- (4.79,5.45);
\draw (4.79,5.45)-- (-0.12,1.88);
\draw (-0.12,1.88)-- (-2,7.66);
\begin{scriptsize}
\draw[color=black] (2.22,7.12) node {$1$};
\draw[color=black] (2.22,4.08) node {$2$};
\draw[color=black] (-0.68,5.02) node {$3$};
\end{scriptsize}
\end{tikzpicture}
\caption{$(m+2)$-angulation $\Delta$}
\label{fig:ex}
\end{figure}
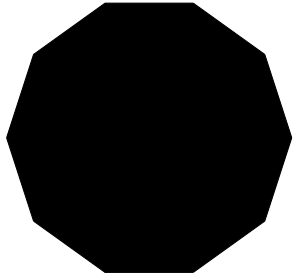

The associated quiver arising from definition \ref{def:grquiv} is the following:

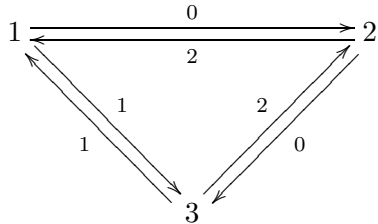
\[\xymatrix{1 \ar@<0.5ex>^0[rrrr] \ar@<0.5ex>^{1}[ddrr] & & & & 2 \ar@<0.5ex>^2[llll] \ar@<0.5ex>^{0}[ddll] \\ & & & & \\ & & 3 \ar@<0.5ex>^{1}[uull] \ar@<0.5ex>^{2}[uurr] & &}\]

We call by $\alpha$ the arrow from $1$ to $2$, $\beta$ the arrow from $2$ to $3$, and $\gamma$ the arrow from $3$ to $1$. Then the associated potential os the sum of all cycles of degree $1$, which means that $W=\gamma \beta \alpha$ with graded differential $\beta \alpha$.

If we flip the $(m+2)$-angulation, we obtain figure \ref{fig:ex2}:

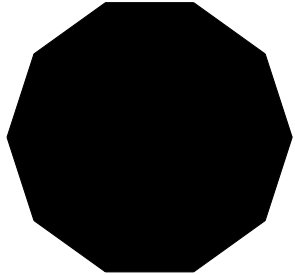
\begin{figure}[!h]
\centering
\begin{tikzpicture}[scale=0.5]
\fill[fill=black,fill opacity=0.1] (-0.12,1.88) -- (2.2,1.88) -- (4.08,3.24) -- (4.79,5.45) -- (4.08,7.66) -- (2.2,9.02) -- (-0.12,9.02) -- (-2,7.66) -- (-2.71,5.45) -- (-2,3.24) -- cycle;
\draw (-0.12,1.88)-- (2.2,1.88);
\draw (2.2,1.88)-- (4.08,3.24);
\draw (4.08,3.24)-- (4.79,5.45);
\draw (4.79,5.45)-- (4.08,7.66);
\draw (4.08,7.66)-- (2.2,9.02);
\draw (2.2,9.02)-- (-0.12,9.02);
\draw (-0.12,9.02)-- (-2,7.66);
\draw (-2,7.66)-- (-2.71,5.45);
\draw (-2.71,5.45)-- (-2,3.24);
\draw (-2,3.24)-- (-0.12,1.88);
\draw (-0.12,9.02)-- (4.79,5.45);
\draw (-0.12,1.88)-- (-2,7.66);
\draw (-2,7.66)-- (4.08,3.24);
\begin{scriptsize}
\draw[color=black] (2.22,7.12) node {$1$};
\draw[color=black] (-0.68,5.02) node {$3$};
\draw[color=black] (1.08,5.88) node {$2$};
\end{scriptsize}
\end{tikzpicture}
\caption{Mutation at vertex $2$ of the initial $(m+2)$-angulation $\Delta$}
\label{fig:ex2}
\end{figure}

which associated quiver is

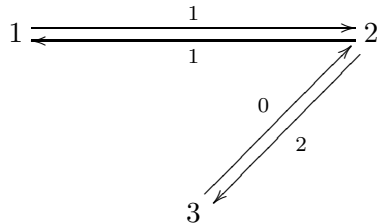
\[\xymatrix{1 \ar@<0.5ex>^1[rrrr] & & & & 2 \ar@<0.5ex>^1[llll] \ar@<0.5ex>^{2}[ddll] \\ & & & & \\ & & 3 \ar@<0.5ex>^{0}[uurr] & &}\]

Now, let us apply Oppermann mutation of the first quiver. We obtain this one:

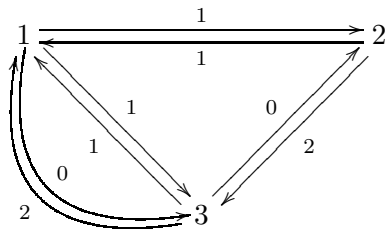
\[\xymatrix{1 \ar@<0.5ex>^1[rrrr] \ar@<0.5ex>^{1}[ddrr] \ar@/_3pc/@<0.5ex>^{0}[ddrr] & & & & 2 \ar@<0.5ex>^1[llll] \ar@<0.5ex>^{2}[ddll]  \\ & & & & \\ & & 3 \ar@<0.5ex>^{1}[uull] \ar@<0.5ex>^{0}[uurr] \ar@/^3pc/@<0.5ex>^{2}[uull] & &}\]

If we factorize through the ideal $(\gamma \alpha^{-1})$, then this quiver is quasi isomorphic to the previous quiver. The new potential is also $0$.

%

\section{Construction of the higher cluster category from a graded quiver with potential}

In this section, we build the higher cluster category arising from a graded quiver with superpotential. We use the results of Keller and Guo established in the preliminaries. In this section, $(Q,W)$ denotes a graded quiver with potential (possibly arising from an $(m+2)$-angulation).

\begin{defi}
From the graded quiver $Q$, we define $\overline{Q}$ as the following quiver. The vertices are the vertices of $Q$, and the arrows are:
\begin{enumerate}
\item If $\alpha$ is an arrow of $Q$ of degree $-r$, then add the same arrow in $\overline{Q}$.
\item If there is $\alpha:\xymatrix@1{i\ar[r]^{-r} & j}$ then add $\alpha^*:\xymatrix@1{j\ar[r]^{-r-m} & i}$ modulo $m$.
\item To each vertex, add a loop of degree $-(m+1)$.
\end{enumerate}
\end{defi}

The graduations graphically correspond to the sides between the edges in an $(m+2)$-angulation.

\begin{defi}\cite[Section 6.2]{Keldefo}
Let $\Gamma^{m+2}(Q,W)$ be the Ginzburg dg category associated with $\overline{Q}$ and $W$. It is endowed with the unique differential which
\begin{enumerate}
\item vanishes on $Q$,
\item takes an element $\alpha^*$ to the cyclic derivative $\partial_{·\alpha^*}W$,
\item takes a loop $t$ at $i$ to \[(-1)^n t (\sum_{\substack{v_i \text{ incident}\\ \text{arrow of }i}} [v_i,v_i^*]) t. \]
\end{enumerate}
\end{defi}

\begin{rmk}
It has been said that we could take the quiver with arrows of graduation $-r$, or the inverse arrows with graduation $-r-m$. Actually, this is the same in the definition of the generalized higher cluster category. Indeed, Keller in \cite{Keldefo} defines the bimodule
\[ Q \bigoplus \check{Q}[m] \bigoplus \mathcal{R}[m+1]. \]
Then if we choose an orientation, some arrows will be in $Q$, and their inverses will lie in $\check{Q}[m]$. If we take the inverse orientation, it does not change the bimodule, then the category.
\end{rmk}

Let us now see the properties of this category.

\begin{theo}\cite{G}
Suppose that $H_0(\Gamma_{m+2}(Q,W))$ is finite-dimensional. Then the generalized $m$-cluster-category
\[ \mathcal{C}^m_{(Q,W)}=\mathrm{per}~\Gamma_{m+2}(Q,W) / \mathcal{D}^b\Gamma_{m+2}(Q,W) \] is Hom-finite and $(m+1)$-Calabi-Yau. Moreover, the image of the free module $\Gamma_{m+2}(Q,W)$ in $\mathcal{C}^m_{(Q,W)}$ is an $m$-cluster-tilting object whose endomorphism algebra is isomorphic to $H_0(\Gamma_{m+2}(Q,W))$.
\end{theo}

In this way, we can define a generalized higher cluster algebra arising from an $(m+2)$-angulation.

\begin{defi}
Let $\Delta$ be an $(m+2)$-angulation arising from an unpunctured surface. Using definition \ref{def:grquiv}, we introduce the graded quiver with superpotential $(Q_\Delta,W_\Delta)$. Then, the generalized higher cluster algebra arising from $\Delta$ is \[ \mathcal{C}^m_{(Q_\Delta,W_\Delta)}=\mathrm{per}~\Gamma_{m+2}(Q_\Delta,W_\Delta) / \mathcal{D}^b\Gamma_{m+2}(Q_\Delta,W_\Delta).\]
\end{defi}

If we draw the information we have in a diagram, it would be the following one:

\[\xymatrix{\Delta \ar[r] \ar[d] & (Q_\Delta,W_\Delta) \ar[r] \ar[d] & \Gamma_{m+2}(Q_\Delta,W_\Delta) \ar[r] \ar[d] & \mathcal{C}^m_{(Q_\Delta,W_\Delta)} \ar[d] \\
\mu_i(\Delta) \ar[r] & Q_{\mu_i(\Delta)}=\mu_i(Q_\Delta) \ar[r] & \Gamma'=\Gamma_{m+2}(\mu_i(Q_\Delta),\mu_i(W_\Delta)) \ar[r] & \mathcal{C}^m_{\Gamma'}}\]

Both squares on the right commute thanks to Oppermann in \cite{O}. The previous section shows the commutativity of the left square.

\begin{rmk}
Unfortunately, it is not possible to associate a generalized higher cluster category to an unpunctured Riemann surface, since all $(m+2)$-angulations are not flip equivalent. A simple counterexample of type $\tilde{A}$ could be the following one:

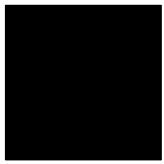
\begin{figure}[!h]
\centering
\begin{tikzpicture}[scale=0.5]
\fill[fill=black,fill opacity=0.1] (-2.26,8.44) -- (-2.26,4.34) -- (1.84,4.34) -- (1.84,8.44) -- cycle;
\fill[fill=black,fill opacity=0.1] (-1.08,7.3) -- (-1.08,5.6) -- (0.62,5.6) -- (0.62,7.3) -- cycle;
\draw (-2.26,8.44)-- (-2.26,4.34);
\draw (-2.26,4.34)-- (1.84,4.34);
\draw (1.84,4.34)-- (1.84,8.44);
\draw (1.84,8.44)-- (-2.26,8.44);
\draw (-1.08,7.3)-- (-1.08,5.6);
\draw (-1.08,5.6)-- (0.62,5.6);
\draw (0.62,5.6)-- (0.62,7.3);
\draw (0.62,7.3)-- (-1.08,7.3);
\draw (-2.26,8.44)-- (-1.08,7.3);
\draw (1.84,8.44)-- (0.62,7.3);
\draw (1.84,4.34)-- (0.62,5.6);
\draw (-2.26,4.34)-- (-1.08,5.6);
\end{tikzpicture}
\end{figure}

This figure is not equivalent by flip to the following one:

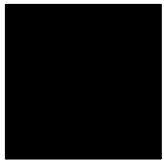
\begin{figure}[!h]

\begin{tikzpicture}[scale=0.5]
\fill[fill=black,fill opacity=0.1] (-2.26,8.44) -- (-2.26,4.34) -- (1.84,4.34) -- (1.84,8.44) -- cycle;
\fill[fill=black,fill opacity=0.1] (-1.08,7.3) -- (-1.08,5.6) -- (0.62,5.6) -- (0.62,7.3) -- cycle;
\draw (-2.26,8.44)-- (-2.26,4.34);
\draw (-2.26,4.34)-- (1.84,4.34);
\draw (1.84,4.34)-- (1.84,8.44);
\draw (1.84,8.44)-- (-2.26,8.44);
\draw (-1.08,7.3)-- (-1.08,5.6);
\draw (-1.08,5.6)-- (0.62,5.6);
\draw (0.62,5.6)-- (0.62,7.3);
\draw (0.62,7.3)-- (-1.08,7.3);
\draw (-2.26,4.34)-- (0.62,5.6);
\draw (1.84,4.34)-- (0.62,7.3);
\draw (1.84,8.44)-- (-1.08,7.3);
\draw (-2.26,8.44)-- (-1.08,5.6);
\end{tikzpicture}
\end{figure}

This is why in his article, Torkildsen tells that we cannot hang an arc from an outside vertex $i$ to an inside vertex $j$ if $i \not\equiv j \mod m$.

Note that in case $m=1$, it is known that any two triangulations are flip equivalent. However, Amiot showed in an appendix of \cite{CS} that, if $\mathcal{C}_\Delta$ (respectively $\mathcal{C}_{\Delta'}$) is a category associated with a triangulation $\Delta$ (respectively $\Delta'$), then the equivalence between $\mathcal{C}_\Delta$ and $\mathcal{C}_{\Delta'}$ is not canonical.
\end{rmk}

\nocite{MP}

\bibliographystyle{alpha}
\bibliography{biblio}

\end{document}